\documentclass[11pt]{amsart}

\title{Good orbital integrals}
\date{10 Sept 2004}

\author{Clifton Cunningham and Thomas C. Hales}

\address{Department of Mathematics, University of Calgary,
Alberta, Canada, T2N 1N4} \email{cunning@math.ucalgary.ca}

\address{Department of Mathematics, University of Pittsburgh,
Pittsburgh, PA 15260} \email{hales@pitt.edu}

\subjclass{22E50, 14F42} \keywords{Orbital integrals, local
constancy, motivic integration, Fundamental Lemma}
\issueinfo{00}
{0}         
{0}         
{}      
\copyrightinfo{2003}
{Cunningham and Hales}

\newtheorem{theorem}{Theorem}[section]
\newtheorem{proposition}[theorem]{Proposition}
\newtheorem{conjecture}[theorem]{Conjecture}
\newtheorem{corollary}[theorem]{Corollary}
\newtheorem{lemma}[theorem]{Lemma}
\theoremstyle{definition}
\newtheorem{definition}[theorem]{Definition}
\newtheorem{example}[theorem]{Example}
\theoremstyle{remark}
\newtheorem{remark}[theorem]{Remark}
\numberwithin{equation}{theorem}

\def\R{{\mathbb{R}}}
\def\C{{\mathbb{C}}}
\def\F{{F}}
\def\A{{\mathfrak{O}_\F}}
\def\pFA{{\mathfrak{p}_\F}}
\def\f{{\mathbb{F}_q}}
\def\barf{{\bar{\mathbb{F}}_q}}

\def\GF{{G}}
\def\AG{{A_G}}
\def\TF{{T}}

\def\LF{{L}}
\def\UF{{U}}
\def\PF{{P}}

\def\lF{{\mathfrak{l}}}
\def\uF{{\mathfrak{u}}}

\def\tF{{\mathfrak{t}}}

\def\gF{{\mathfrak{g}}}

\def\t{{\mathfrak{t}}}
\newcommand{\ellip}{{e}}

\def\orbint{{\phi}}

\newcommand{\Mf}{{\bar\GF}}
\newcommand{\mf}{{\bar\gF}}
\def\vphi{{\varphi}}
\def\hphi{{\hat{\varphi}}}
\def\orbit{{\mathcal{O}}}
\def\X{{\mathcal{X}}}
\def\Y{{\mathcal{Y}}}
\def\Z{{\mathcal{Z}}}
\def\calL{{\mathcal{L}}}
\def\calC{{\mathcal{C}}}
\def\calA{{\mathcal{A}}}
\def\calB{{\mathcal{B}}}
\def\calF{{\mathcal{F}}}

\def\Kill#1#2{{\langle{#1},{#2}\rangle}}
\def\kill#1#2{{\Lambda\left({#1},{#2}\right)}}
\def\i{{i}}

\def\B{{\mathcal{B}}}

\def\Fourier{{\mathcal{F}}}
\def\Sch#1{{C^\infty_c(#1)}}
\def\Ad{{\rm Ad}}
\def\rightcoset{{ / }}
\def\leftcoset{{\backslash}}
\def\vol{{\operatorname{vol}}}
\def\meas{{\operatorname{mes}}}
\def\tq{{\ \vert\ }}
\def\ad{{\rm ad}}
\def\nilp{{nil}}
\def\Lie{{\operatorname{Lie}}}
\def\iso{{\ \cong\ }}

\def\abs#1{{\vert #1\vert}}

\newcommand{\G}{{\mathbb{G}}}

\newcommand{\T}{{\mathbb{T}}}
\def\bb#1{{\mathbb #1}}
\newcommand{\ring}[1]{\mathbb{#1}}
\def\Fq{{\mathbb F}_q}
\def\op#1{\operatorname{#1}}
\def\lie#1{{\mathfrak{#1}}}
\def\liegood#1{{\lie#1}}
\def\orb#1#2{{\phi_{\lie{#1},\lie{#2}}}}
\def\mot{{K_0^{mot}(\op{Var}_k)[\ring{L}^{-1}]\otimes\ring{Q}}}
\def\vol{\operatorname{vol}}
\def\ac{\operatorname{ac}}
\def\sign{\operatorname{sign}}

\def\Ad{\operatorname{Ad}}
\def\ad{\operatorname{ad}}

\def\chr{\operatorname{char}}
\def\rank{\operatorname{rank}}
\def\rq{\text{'}}
\def\lb{\text{\#}}

\begin{document}

\newcommand{\FF}{{\mathbb{F}}}

\begin{abstract}
This paper concerns a class of orbital integrals in Lie algebras
over $p$-adic fields.  The values of these orbital integrals at
the unit element in the Hecke algebra count points on varieties
over finite fields.  The construction, which is based on motivic
integration, works both in characteristic zero and in positive
characteristic. As an application, the Fundamental Lemma for this
class of integrals is lifted from positive characteristic to
characteristic zero.  The results are based on a formula for
orbital integrals as distributions inflated from orbits in the
quotient spaces of the Moy-Prasad filtrations of the Lie algebra.
This formula is established by Fourier analysis on these quotient
spaces.
\end{abstract}

\maketitle

\def\today{\ifcase\month\or
    January\or February\or March\or April\or May\or June\or
    July\or August\or September\or October\or November\or December\fi
    \space\number\day, \number\year}

\maketitle

\section*{Introduction}\label{section: introduction}

It has been clear to researchers for many years that orbital
integrals on $p$-adic groups are geometrical objects.  However, it
has taken many years to make this observation precise.  When the
local field has positive characteristic, Kottwitz, Goresky, and
MacPherson give a geometrical description of the orbital integrals
of the unit element in the Hecke algebra \cite{GKM-fl}.

Another approach to this problem is suggested by motivic
integration.  Motivic integration may be viewed as a
geometrization of ordinary $p$-adic integration.  This is the path
followed by this paper.

One advantage of this approach is that it works equally well in
all characteristics.  This allows us to lift the beautiful recent
work of Goresky, Kottwitz, Laumon, and MacPherson to
characteristic zero -- at least for the special class of
semi-simple elements that we consider.

A limitation of motivic integration is that the domain of
integration is restricted to a special algebra of sets, called
definable sets (in the sense of first order logic). The main
problem we face is that $p$-adic orbits are not definable sets.
However, orbital integrals are locally constant functions. This
allows us to replace each orbital integral by an average of
orbital integrals over a neighborhood.

The key question is then whether orbital integrals are constant on
definable sets. This turns out to be the case, at least for the
class of orbital integrals that we study in this paper. In fact,
the relation between definable sets and the local constancy of
orbital integrals is rather striking.  What we find is that is the
largest neighborhood of an element on which (we are able to show
that) the orbital integrals are constant coincides precisely with
the smallest neighborhood that is definable.  In brief, we find
that motivic integration is perfectly adapted to the study of
orbital integrals.

The first section of this paper uses local Fourier analysis to
determine a neighborhood of certain semi-simple elements on which
orbital integrals are constant. The main results on the local
constancy of orbital integrals are Theorem~\ref{theorem: main} and
Corollary~\ref{cor:unit element}. The remainder of the paper shows
that this neighborhood is definable and then applies the machinery
of motivic integration to give a geometric interpretation of
orbital integrals. This geometric interpretation of orbital
integrals, Theorem~\ref{thm:variety}, is the main result of the
paper.

This paper carries this project through for a significant special
case, although we assume the residual characteristic of the
$p$-adic field is large. In Section \ref{section: formula} the
local Fourier analysis is treated for all connected reductive
groups, but with restrictions on the valuations of the roots of
the semi-simple elements. (These are the {\it good\/} elements
that appear in the title.) The remainder of this paper restricts
further to classical groups, and places some further restrictions
on the valuations of the roots.

As an application of Theorem~\ref{thm:variety}, we observe in
Corollary~\ref{cor:F} that the Fundamental Lemma for this class of
integrals is lifted from positive characteristic to characteristic
zero.

It is our expectation that the results should generalize to all
reductive groups and all semi-simple elements without restriction
on the valuations of the roots; however, this has not yet been
carried out.

Ju-Lee Kim has pointed out that the explicit local constancy of
orbital integrals on the set of good semisimple elements can be
deduced directly from the paper \cite{KM}.  An appendix by Ju-Lee
Kim indicates how to do this.  Ju-Lee Kim has recently proved more
general results about the local constancy of orbital integrals
(\cite[Theorem~9.2.2]{KM2}).

After this paper was already submitted for publication, we
received the preprint \cite{WC}, which treats a similar topic as
our paper. Corollary~\ref{cor:F} is a special case of
\cite[thm7.2]{WC}.

We are happy to acknowledge the assistance of Jeffrey Adler and to
thank him for many helpful conversations. We would also like to
thank Fran\c{c}ois Loeser for his support of this project.

\tableofcontents

\section{A Formula for Good Orbital Integrals}\label{section:
formula}

\subsection{Preliminary remarks}\label{section:
preliminaries}

Let $\F$ be a $p$-adic field with ring of integers $\A$, prime
ideal $\pFA$ and residue field $\Fq$.  Throughout the paper, there
will be certain mild restrictions on $p$, the characteristic of
the residue field of $F$; these conditions are met my assuming
that $p$ is `sufficiently large'.

Let $|x| = |x|_F$ be the normalized absolute value on $F$. (We
also use $|C|$ for the cardinality of a set $C$.  If
$\xi=(\xi_1,\ldots,\xi_m)$ is an $m$-tuple of variables, then we
write $|\xi|$ for the length $m$ of the tuple.  The context will
make it clear which is intended.)

Let $\GF$ be the group of $\F$-rational points on a connected
reductive algebraic group $\G$ defined over $\F$. Integration on
$\GF$ and its measurable subsets will always be taken with respect
to the same Haar measure and the notation `$\meas$' will refer to
that measure. Let $\gF$ denote Lie algebra of $\GF$. Every
integral over $\gF$ and its measurable subsets will be taken with
respect to the same Haar measure and the notation `$\vol$' will
refer to that measure.

The Bruhat-Tits building for $\GF$ will be denoted $\B(\GF)$; we
refer the reader to \cite{Bruhat-Tits II} for the definition.  A
torus $\T\subseteq\G$ defined over $\F$ is \emph{tamely ramified}
if the splitting field for $\T$ over $\F$ has ramification index
prime to $p$. In this case, using results of \cite{L}, we choose
an embedding of the Bruhat-Tits building $\B(\TF)$ for the group
$\TF$ of $\F$-rational points on $\T$ into the Bruhat-Tits
building $\B(\GF)$ for $\GF$ by way of a toral map; as we use only
the image of this embedding, which does not depend on the choice
just made, all results in this paper are independent of the
choice. We will also view the building $\B(\LF)$ for the group of
$\F$-rational points on a Levi subgroup $\mathbb{L}\subset \G$
defined over $\F$ as a simplicial subcomplex of $\B(\GF)$.

The reader is referred to \cite{Moy-Prasad} for the definition of
the subgroups $\GF_{x,r}$  and $\GF_{x,r^+}$ of $\GF$ (resp. the
lattices $\gF_{x,r}$ and $\gF_{x,r^+}$ of $\gF$) where $x$ is any
point in $\B(\GF)$ and $r$ is any non-negative real number (resp.
any real number). We also write $\gF_{r}$ for the union of the
spaces $\gF_{x,r}$, as $x$ ranges over all $\B(\GF)$.

Let $\Sch{\gF}$ denote the convolution $\C$-algebra of locally
constant complex-valued functions on $\gF$ with compact support.
For any lattice $\mathcal{L}$ in $\gF$, we write
$\Sch{\mathcal{L}}$ for the space of functions in $\Sch{\gF}$
supported by $\mathcal{L}$; if $\mathcal{L}^\prime$ is a
sublattice in $\mathcal{L}$, then
$\Sch{\mathcal{L}/\mathcal{L}^\prime}$ will denote the vector
space of elements of $\Sch{\gF}$ supported by $\mathcal{L}$ which
are constant on the $\mathcal{L}^\prime$-cosets in $\mathcal{L}$.

For each $x\in \B(\GF)$, we let $\GF_x$ be the parahoric subgroup
associated to $x$ by \cite{Bruhat-Tits II} and let $\Mf_x$ denote
the quotient group $\GF_x\rightcoset\GF_{x,0^+}$. We remark that
$\Mf_x$ is the set of $\Fq$-points of a reductive linear algebraic
group. For each real number $r$ (and any $x\in \B(\GF)$) let
$\mf_{x,r}$ denote $\gF_{x,r}\rightcoset\gF_{x,r^+}$ and let
$\rho_{x,r} : \gF_{x,r} \to \mf_{x,r}$ be the projection map. We
remark that the adjoint action of $\GF$ on $\gF$ restricts to an
action of $\GF_x$ on $\gF_{x,r}$ which in turn induces an action
of $\Mf_x$ on $\mf_{x,r}$.

Let $d_x(X)$ denote the supremum of the set of all $r\in \R$ such
that $X \in \gF_{x,r}$. This defines a function $d_x: \gF \to \R
\cup \left\{+\infty\right\}$ which we refer to as the {\em depth
function at $x$}. The {\em depth} $d(X)$ of $X$ in $\gF$ is the
supremum of the $d_x(X)$ as $x$ ranges over the building
$\B(\GF)$. The depth of a non-zero non-nilpotent element is always
a rational number. The depth of $X$ is infinite exactly when $X$
is nilpotent.

We will make extensive use of the notion of {\em good} elements in
$\gF$, as introduced in \cite[2.2.4]{Adler}.  Accordingly, we
review that definition here. First, let $\T\subseteq\G$ be a
tamely ramified torus defined over $\F$ and let $E$ be a splitting
field for $\T$. Let $\TF(E)$ denote the group of $E$-rational
points on $\T$ and let $\tF\otimes E$ denote the Lie algebra of
$\T(E)$. Note that $\tF\otimes E$ is a split Lie algebra since it
coincides with the Lie algebra of the group of $E$-rational points
on the connected reductive algebraic group $\T
\times_{\operatorname{Spec}(F)} \operatorname{Spec}(E)$, which is
split over $E$ by construction. For any $Y \in \tF\otimes E$ and
any $y\in \B(\T(E))$ we define the depth of $Y$ at $y$ as above,
and likewise define the depth of $Y$ in $\tF\otimes E$ as above.

\begin{definition}\label{definition: good}
An element $X$ of $\gF$ is {\em good} if $X$ is semi-simple, $X$
is contained in a Cartan subalgebra $\tF = \Lie(\TF)$ which is
tamely ramified with splitting field $E/F$, and for every root
$\alpha$ of $\gF$ relative to $\tF$, either $\alpha(X)$ is zero or
the $E$-normalized valuation of $\alpha(X)$ equals the depth of
$X$ in $\tF\otimes E$.
\end{definition}

\begin{remark}\label{remark: good}
It should be noted that the parameterization of the filtrations in
\cite{Adler} differs by a scalar multiple from that of
\cite{Moy-Prasad}. In Definition~\ref{definition: good}, we use
the depth only on a split Lie algebra $\tF\otimes E$; and here the
two parameterizations of the filtrations coincide. In this paper,
we use the parameterization defined in \cite{Moy-Prasad}, so all
results culled from \cite{Adler} are translated accordingly.
\end{remark}

The existence of good elements, assuming $p$ sufficiently large,
is established in \cite{A}.


\subsection{The Fourier transform}\label{section: Fourier transform}

This section fixes a Fourier transform $\Fourier_\gF$ on the
$p$-adic Lie algebra and recalls some well-known elementary
properties of $\Fourier_\gF$. Throughout, $x$ is an arbitrary
element of $\B(\GF)$ and $r$ is an arbitrary real number.

We fix a Killing form $\Kill{\cdot}{\cdot}:\gF\times\gF\to\F$ for
$\gF$ and a non-trivial additive character
$\lambda:\F\to\C^\times$ with conductor $\A$ (that is, $\lambda$
is trivial on $\pFA$ but not trivial on $\A$).  Let $\kill{X}{Y}$
denote the image of $\Kill{X}{Y}$ under $\lambda$.

\begin{lemma}\label{lemma: FT}
Let $1_{Z+\gF_{x,r}} : \gF \to \C$ be the characteristic function
of $Z+\gF_{x,r}$. For any $x\in \B(\GF)$, $r\in \R$ and $Z\in
\gF$,
\begin{equation}\label{equation: FT}
\forall X\in \gF, \quad \int_{\gF} \kill{X}{Y}\ 1_{Z+\gF_{x,r}}
(Y)\, dY
 = \kill{X}{Z} \vol(\gF_{x,r}) 1_{\gF_{x,(-r)^+}}(X).
\end{equation}
\end{lemma}

\begin{proof}
\begin{eqnarray*}
\int_{\gF} \kill{X}{Y}\ 1_{Z+\gF_{x,r}} (Y)\, dY
 &=& \int_{Z+\gF_{x,r}} \kill{X}{Y}\, dY\\
 &=& \int_{\gF_{x,r}} \kill{X}{Z+Y}\, dY\\
 &=& \kill{X}{Z}\int_{\gF_{x,r}} \kill{X}{Y}\, dY.
\end{eqnarray*}
For $p$ sufficiently large (recall this assumption from section
\ref{section: preliminaries}) it follows from the definition of
the lattice $\gF_{x,r}$ and the fact that $\lambda$ has conductor
$\A$ that the set of all $X\in \gF$ for which $\kill{X}{Y}=1$ for
all $Y\in \gF_{x,r}$ is exactly $\gF_{x,(-r)^+}$ (ref:
\cite[\S~4]{Adler}). Thus, when $X$ is an element of
$\gF_{x,(-r)^+}$ we have
\begin{eqnarray*}
\kill{X}{Z}\int_{\gF_{x,r}} \kill{X}{Y}\, dY
 &=& \kill{X}{Z} \vol(\gF_{x,r});
\end{eqnarray*}
on the other hand, if $X$ is not contained in $\gF_{x,(-r)^+}$
then the function $Y \mapsto \kill{X}{Y}$ is non-trivial on the
lattice $\gF_{x,r}$, so
\begin{eqnarray*}
\kill{X}{Z}\int_{\gF_{x,r}} \kill{X}{Y}\, dY
 &=& 0.
\end{eqnarray*}
Therefore,
\begin{equation}
\int_{\gF} \kill{X}{Y}\ 1_{Z+\gF_{x,r}} (Y)\, dY
 = \kill{X}{Z} \vol(\gF_{x,r}) 1_{\gF_{x,(-r)^+}}(X),
\end{equation}
as claimed.
\end{proof}

\begin{proposition}\label{proposition: FT}
If $f: \gF \to \C$ is locally constant with compact support, then
the function $X \mapsto \int_{\gF} \kill{X}{Y}\ f(Y)\, dY$ is also
locally constant with compact support.
\end{proposition}

\begin{proof}
By hypothesis, $f\in \Sch{\gF}$ (as defined in
Section~\ref{section: preliminaries}). Any element of $\Sch{\gF}$
may be expressed as a finite linear combination of functions of
the form $1_{Z+\gF_{x,r}}$. Thus, it is sufficient to prove the
proposition when $f = 1_{Z+\gF_{x,r}}$. From Lemma~\ref{lemma: FT}
we see that $X \mapsto \int_{\gF} \kill{X}{Y}\ 1_{Z+\gF_{x,r}}\,
dY$ is locally constant and compactly supported, and therefore an
element of $\Sch{\gF}$.
\end{proof}

\begin{definition}\label{definition: Fourier transform}
Define $\Fourier_\gF : \Sch{\gF} \to \Sch{\gF}$ by
\begin{equation}
(\Fourier_\gF f)(X) = \int_{\gF} \kill{Y}{X}\ f(Y)\, dY.
\end{equation}
We refer to $\Fourier_\gF$ as the {\em Fourier transform on
$\gF$}. When there is no ambiguity, we write $\widehat{f}$ for
$\Fourier_\gF f$.
\end{definition}

\begin{corollary}\label{corollary: FT bijection}
For any $x\in \B(\GF)$ and for any $r,s\in \R$ with $s\leq r$, the
Fourier transform induces an isomorphism of vector spaces
\begin{equation}
\Fourier_\gF : \Sch{\gF_{x,s}\rightcoset\gF_{x,r^+}}\to
\Sch{\gF_{x,-r}\rightcoset\gF_{x,(-s)^+}}.
\end{equation}
\end{corollary}

\begin{proof}
Any element of $\Sch{\gF_{x,s}\rightcoset\gF_{x,r^+}}$ is a finite
linear combination of functions of the form $1_{Z+\gF_{x,r^+}}$
with $d_x(Z) = s$. Let $r'$ be the largest real number such that
$\gF_{x,r'} = \gF_{x,r^+}$. (Thus, $r'$ is the first `jump point'
greater than $r$.) We have seen in the proof of Proposition
\ref{proposition: FT} (cf: Equation \ref{equation: FT}) that
\[
(\Fourier_\gF 1_{Z+\gF_{x,r'}})(X) = \kill{X}{Z} \vol(\gF_{x,r'})
1_{\gF_{x,(-r')^+}}(X).
\]
Since $\gF_{x,(-r')^+} = \gF_{x,-r}$, we have
\[
(\Fourier_\gF 1_{Z+\gF_{x,r^+}})(X) = \kill{X}{Z}
\vol(\gF_{x,r^+}) 1_{\gF_{x,-r}}(X),
\]
so $\Fourier_\gF$ maps $\Sch{\gF_{x,s}\rightcoset\gF_{x,r^+}}$
into $\Sch{\gF_{x,-r}}$. Suppose now that $X\in \gF_{x,-r}$. If
$X' \in \gF_{x,(-s)^+}$ then $X'\in \gF_{x,-r}$ so
\begin{eqnarray*}
(\Fourier_\gF 1_{Z+\gF_{x,r+}})(X+X')
 &=& \kill{X+X'}{Z}\vol(\gF_{x,r^+}) 1_{\gF_{x,-r}}(X+X')\\
 &=& \kill{X}{Z}\kill{X'}{Z}\vol(\gF_{x,r^+})\\
 &=& \kill{X}{Z}\vol(\gF_{x,r^+}).
\end{eqnarray*}
Thus, $\Fourier_\gF$ maps $\Sch{\gF_{x,s}\rightcoset\gF_{x,r^+}}$
into $\Sch{\gF_{x,-r}/\gF_{x,(-s)^+}}$. Applying the same ideas we
see that $\Fourier_\gF\circ \Fourier_\gF$ maps
$\Sch{\gF_{x,s}\rightcoset\gF_{x,r^+}}$ into
$\Sch{\gF_{x,s}\rightcoset\gF_{x,r^+}}$. To see that $\Fourier_\gF$
is an isomorphism, let $f = 1_{Z+\gF_{x,r^+}}$ and use Equation
\ref{equation: FT} again to see that $(\Fourier_\gF(\Fourier_\gF
f))(X)$ is a scalar multiple (independent of $X$) of $f(-X)$. Since
functions of this form give a basis for
$\Sch{\gF_{x,s}\rightcoset\gF_{x,r^+}}$, the corollary is proved.
\end{proof}


\subsection{Support and mesh}

Recall from Section \ref{section: preliminaries} that we write
$\gF_{r}$ for the union of the spaces $\gF_{x,r}$, as $x$ ranges
over all $\B(\GF)$. Notice that $r\leq r'$ implies $\gF_r
\supseteq \gF_{r'}$.

\begin{definition}\label{definition: depth r Hecke} For any pair
of real numbers $s\leq r$ let $\Sch{\gF}^s_r$ denote the space of
$f\in\Sch{\gF}$ such that the support of $f$ is contained in
$\gF_s$ and the support of $\widehat{f}$ is contained in
$\gF_{-r}$.
We write $\Sch{\gF}_r$ for the union of the $\Sch{\gF}_r^s$ with $s\leq r$.
\end{definition}

\begin{lemma}\label{lemma: depth r Hecke}
Fix $r\in \R$. A compactly supported function $f\in \Sch{\gF}$ is
contained in $\Sch{\gF}^s_r$ if and only if there is a finite set
$\{y_i \tq i\in I\} \subset \B(\GF)$ such that
\begin{equation}
f = \sum_{i\in I} f_i, \qquad \text{ with } f_i \in
\Sch{\gF_{y_i,s} \rightcoset \gF_{y_i,r^+}}. \label{eqn:f
partition}
\end{equation}
\end{lemma}

\begin{proof}
Suppose $f = \sum_{i\in I} f_i$ for some finite set $I$, where
$f_i$ is contained in $\Sch{\gF_{y_i,s} \rightcoset
\gF_{y_i,r^+}}$ for some $y_i\in \B(\GF)$. Then, the support of
$f$ is contained in the union of the $\gF_{y_i,s}$, which is
contained in $\gF_s$. Consider $\widehat{f} = \sum_{i\in I}
\widehat{f_i}$. Applying Corollary~\ref{corollary: FT bijection},
it follows that $\widehat{f_i}$ is an element of
$\Sch{\gF_{y_i,-r}\rightcoset\gF_{y,(-s)^+}}$. In particular, it
follows that the support of $\widehat{f_i}$ is contained in
$\gF_{y_i,-r} \subset \gF_{-r}$, so the support of $\widehat{f}$
is contained in the union of the $\gF_{y_i,-r}$, which is
contained in $\gF_{-r}$. It follows that $f\in \Sch{\gF}^s_r$.

Conversely, fix $f\in \Sch{\gF}^s_r$. Since the support of
$\widehat{f}$ is compact, there is a finite set $\{ y_i \tq i \in
I \}$ such that the support of $\widehat{f}$ is covered by $\{
\gF_{y_i,-r} \tq i \in I\}$. Let $\mu = \sum_{i\in I} \mu_i$ be a
partition of unity for $\cup_{i\in I} \gF_{y_i,-r}$; thus, in
particular, $\mu_i \in \Sch{\gF_{y_i, -r}}$ for each $i\in I$, and
if $j \ne i$ then $\mu_j$ and $\mu_i$ are supported by disjoint
sets. Let $f_i = f \mu_i$ and observe that $f_i \in
\Sch{\gF_{y_i,-r}}$. Then $\widehat{f} = \sum_{i\in I}
\widehat{f_i}$. Since $\widehat{f_i}$ is locally constant,
$\widehat{f_i}$ is an element of
$\Sch{\gF_{y_i,-r}/\gF_{y_i,(-s_i)^+}}$, for some $s_i \in \R$.
Using Corollary~\ref{corollary: FT bijection} again, it follows
that $f_i$ is an element of $\Sch{\gF_{y_i,s_i} \rightcoset
\gF_{y_i,r^+}}$. Since the support of $f$ (and therefore $f_i$) is
contained in $\gF_s$ by hypothesis, we have $s_i \leq s$ for each
$i\in I$. Since $\Sch{\gF_{y_i,-r}/\gF_{y_i,(-s_i)^+}} \subseteq
\Sch{\gF_{y_i,-r}/\gF_{y_i,(-s)^+}}$, it follows that $f_i \in
\Sch{\gF_{y_i,-r}/\gF_{y_i,(-s)^+}}$, for each $i\in I$.
\end{proof}

As an immediate consequence we have the following.
\begin{proposition}\label{proposition: depth r Hecke}
For all $s\leq r$, $\Sch{\gF}^s_r$ is a vector space and the
Fourier transform $\Fourier_\gF$ induces an automorphism
\begin{equation}
\Fourier_\gF : \Sch{\gF}^s_r \to \Sch{\gF}^{-r}_{-s}.
\end{equation}
\end{proposition}

Notice also that $\Sch{\gF}^{s'}_r \supseteq
\Sch{\gF}^s_r$ when $s'\leq s \leq r$ then and $\Sch{\gF}^s_{r}
\subseteq \Sch{\gF}^s_{r'}$ when $s \leq r \leq r'$.
Thus, $\Sch{\gF}_r \subseteq \Sch{\gF}_{r'}$ when $r\leq r'$.

\begin{remark}\label{remark: unit}
Note that $1_{\gF(\A)}$ is an element $\Sch{\gF}_0^0$ and therefore an element of $\Sch{\gF}_0$.
It follows that $1_{\gF(\A)}$ is an element of $\Sch{\gF}_r$ for every $r\geq 0$.
\end{remark}


\subsection{Relative Fourier transform}\label{section: relative
Fourier transform}

This section considers a Fourier transform on the space of
complex-valued functions on $\mf_{x,r}$. The Fourier transform
will then be related to the Fourier transform on the $p$-adic Lie
algebra by inflation. Throughout this section, $x$ is an arbitrary
element of $\B(\GF)$ and $r$ is an real number.

Define $\Lambda_{x,r}:\mf_{x,r}\times\mf_{x,-r}\to\C$ by
\begin{equation}
\Lambda_{x,r}(\X,\Y) = \kill{X}{Y},
\end{equation}
where $X$ is any representative for $\X$ and $Y$ is any
representative for $\Y$. If $X_1$ and $X_2$ are elements of
$\gF_{x,r}$ with $X_1-X_2\in \gF_{x,r^+}$ and if $Y_1$ and $Y_2$
are elements of $\gF_{x,-r}$ with $Y_1-Y_2\in \gF_{x,(-r)^+}$,
then $\kill{X_1}{Y_1}$ equals $\kill{X_2}{Y_2}$, as we have seen
in the proof of Proposition~\ref{proposition: FT}. Thus,
$\Lambda_{x,r}$ is well-defined. The function $\Lambda_{x,r}$
defines a perfect pairing between $\mf_{x,r}$ and $\mf_{x,-r}$.
Note also that $\Lambda_{x,-r}(\Y,\X) = \Lambda_{x,r}(\X,\Y)$.

\begin{definition}\label{definition: finite Fourier transform}
Let $\C(\mf_{x,r})$ denote the space of complex-valued functions
on $\mf_{x,r}$. Define $\Fourier_{x,r} : \C(\mf_{x,r}) \to
\C(\mf_{x,-r})$ by
\begin{equation}
(\Fourier_{x,r}\, \vphi)(\X) = \sum_{\Y \in \mf_{x,r}}
\Lambda_{x,r}(\Y,\X) \vphi(\Y).
\end{equation}
We refer to $\Fourier_{x,r}$ as the {\em finite Fourier transform}
for the pair $(x,r)$. When there is no ambiguity, we write $\hphi$
for $\Fourier_{x,r}\, \vphi$.
\end{definition}

An elementary calculations shows that, for any
$\vphi\in\C(\mf_{x,r})$,
\begin{equation}
    \forall
    \X\in\mf_{x,r},\qquad (\Fourier_{x,-r}(
    \Fourier_{x,r} \vphi)) (\X) = {\vert\mf_{x,r}\vert}\
    \vphi(-\X).
\end{equation}
It is therefore common to define the finite Fourier transform by
introducing the factor ${\vert\mf_{x,r}\vert}^{-1/2}$; as this
does not simplify our main result Theorem~\ref{theorem: main}, we
have not followed that convention here.

\begin{definition}\label{definition: inflation}
For each $\vphi$ in $\C(\mf_{x,r})$, define $\vphi_{x,r} \in \Sch{\gF}$ by
\begin{equation}
\vphi_{x,r}(Y) =
\begin{cases}
(\vphi \circ \rho_{x,r})(Y),& \forall Y \in \gF_{x,r}\\
0,& \forall Y \not\in \gF_{x,r}.
\end{cases}
\end{equation}
The map taking $\vphi$ to $\vphi_{x,r}$ is commonly referred to as
{\em inflation} from $\C(\mf_{x,r})$ to $\Sch{\gF}$.
\end{definition}

\begin{proposition}\label{proposition: inflation and FT}
For any $\vphi \in \C(\mf_{x,r})$,
\begin{equation}
\widehat{\vphi_{x,r}} = \vol(\gF_{x,r^+})\ \hphi_{x,-r}.
\end{equation}
\end{proposition}

\begin{proof}
We will show $\widehat{\vphi_{x,r}}(Y) = \vol(\gF_{x,r^+})\
\hphi_{x,-r}(Y)$ for each $Y\in\gF$ by considering two cases.
First, suppose $Y \in \gF_{x,-r}$. Then,
\begin{equation}
    \widehat{\vphi_{x,r}}(Y) = \int_{\gF_{x,r}} \kill{Z}{Y}\
    \vphi_{x,r}(Z)\, dZ.
    \label{eqn:LpR}
\end{equation}
Recall that
\begin{equation}
\forall
    Z\in\gF_{x,r},\qquad
\kill{Z}{Y} = \Lambda_{x,r}(\rho_{x,r}(Z),\rho_{x,-r}(Y)).
\end{equation}
We now pick a set $\tilde{\gF}_{x,r}$ of representatives for
$\mf_{x,r}$ and write $Z = Z_{r} + Z^\prime$ where $Z_{r} \in
\tilde{\gF}_{x,r}$ and $Z^\prime \in \gF_{x,r^+}$.  Then
$\rho_{x,r}(Z) = \rho_{x,r}(Z_r)$ and $\kill{Z}{Y} =
\kill{Z_r}{Y}\ \kill{Z^\prime}{Y}$. Thus,
\begin{equation}
\kill{Z}{Y} = \Lambda_{x,r}(\rho_{x,r}(Z_r),\rho_{x,-r}(Y)).
\end{equation}
Combining this with Equation~\ref{eqn:LpR} we have
\begin{eqnarray*}
\widehat{\vphi_{x,r}}(Y)
 &=& \int_{\gF_{x,r}} \kill{Z}{Y}\ \vphi_{x,r}(Z)\, dZ\\
 &=& \int\limits_{\gF_{x,r} / \gF_{x,r^+}}\ \int\limits_{\gF_{x,r^+}}
\Lambda_{x,r}(\rho_{x,r}(Z_r),\rho_{x,-r}(Y))\
\vphi(\rho_{x,r}(Z_r))\, dZ^\prime\, dZ_r,
\end{eqnarray*}
where $dZ_r$ denotes the quotient measure on $ \gF_{x,r} /
\gF_{x,r^+}$. Notice that this integrand does not depend on
$Z^\prime$ and that $\int_{\gF_{x,r^+}}dZ^\prime =
\vol(\gF_{x,r^+})$. Thus,
\begin{eqnarray*}
&& \int\limits_{\gF_{x,r} / \gF_{x,r^+}}\
\int\limits_{\gF_{x,r^+}}
\Lambda_{x,r}(\rho_{x,r}(Z_r),\rho_{x,-r}(Y))\
\vphi(\rho_{x,r}(Z_r))\, dZ^\prime\, dZ_r\\
 &=& \vol(\gF_{x,r^+})
\sum_{\Z \in \mf_{x,r}}
        \Lambda_{x,r}(\Z,\rho_{x,-r}(Y))\ \vphi(\Z)\\
 &=& \vol(\gF_{x,r^+})\ \hphi_{x,-r}(Y).
\end{eqnarray*}
This proves the proposition in the first case.

Next, suppose $Y \not\in \gF_{x,-r}$ and notice that it follows
that $\hphi_{x,-r}(Y) = 0$. Let $s = - d_x(Y)$ and note that $s >
r$. As above, pick a set of representatives for
$\gF_{x,r}/\gF_{x,s}$ and, for each $Z\in \gF_{x,r}$, write $Z =
Z_r + Z_s$, where $Z_r$ is from that set of representatives and
$Z_s \in \gF_{x,s}$. Then,
\begin{equation}
\kill{Z}{Y}
=\Lambda_{x,s}\left(\rho_{x,s}(Z_s),\rho_{x,-s}(Y)\right)\
\kill{Z_r}{Y}.
\end{equation}
Thus,
\begin{eqnarray*}
\widehat{\vphi_{x,r}}(Y)
 &=& \int_{\gF_{x,r} / \gF_{x,s} } \kill{X}{Y_r}\
\vphi(\rho_{x,r}(Z_r))\, dZ_r \\
&& \times
\int_{\gF_{x,s}}
\Lambda_{x,s}(\rho_{x,s}(Z_s),\rho_{x,-s}(Y))\, dZ_s.
\end{eqnarray*}
Since $\rho_{x,-s}(Y)\not= 0$, the second integral is $0$ and
therefore $\widehat{\vphi_{x,r}}(Y)= 0$. This proves the
proposition in the second case, and therefore finishes the proof
of Proposition~\ref{proposition: inflation and FT}.
\end{proof}

Proposition~\ref{proposition: inflation and FT} states that the
Fourier transform commutes with inflation, up to a multiple.


\subsection{Gauss\ integrals}\label{section: gauss integrals}

This section introduces our main technique for the study of
regular semi-simple orbital integrals.

\begin{definition}\label{definition: Gauss integral}
For any point $x$ in the building for $\GF$, define $\i_x : \gF
\times \gF \to \C$ by
\begin{equation}
\i_x(X,Y) = \int_{\GF_x} \kill{\Ad(g)X}{Y}\, dg.
\end{equation}
We refer to $i_x(X,Y)$ as a \emph{Gauss\ integral}. We will
sometimes write $\i_{x,X}$ for the function $Y \mapsto i_x(X,Y)$.
\end{definition}

We begin with three lemmas due to Jeffrey Adler.  We write $G_X$
for the centralizer of $X$ in $G$.

\begin{lemma}\label{lemma: Adler 1}
If $X$ is a good regular element of $\gF$, and if $x$ is in the
building for $\GF_X$ in $\GF$, then
\begin{equation}
    \forall Y \in \t^\perp,\qquad d_x([X,Y]) = d_x(X) + d_x(Y),
\end{equation}
where $\t^\perp$ denotes the subspace of $\gF$ perpendicular to
$\t = \Lie\,\GF_X$ with respect to the killing form on $\gF$.
\end{lemma}

\begin{proof}
The result follows from \cite[2.3.1]{Adler}.
\end{proof}

\begin{lemma}\label{lemma: Adler 2}
Suppose $x$ is a point in the building for a maximal torus $\TF$
in $\GF$. For each $t > 0$ there is a diffeomorphism $e_{x,t} :
\gF_{x,t} \to \GF_{x,t}$ such that if $Z \in \gF_{x,t}$, then
\begin{equation}
    \forall Y\in\gF,\qquad\Ad(e_{x,t}(Z)) Y \in Y + \ad(Z) Y +
    \gF_{x,2t+d_x(Y)}.
\end{equation}
\end{lemma}

\begin{proof}
This homeomorphism is defined at the end of \cite[\S~1.5]{Adler}.
The properties above follow from \cite[1.6.3]{Adler} and
\cite[1.6.4]{Adler}. We have also used \cite[1.6.6]{Adler}.
\end{proof}

\begin{lemma}\label{lemma: Adler 3}
Suppose that $X$ is tamely ramified and that $x$ is a point in the building
for
$\GF_X$ in $\GF$. Let $r = d_x(X)$. Then $X + \gF_{x,r^+}$ contains no
nilpotent
elements.
\end{lemma}

\begin{proof}
This is \cite[1.9.5]{Adler}.
\end{proof}

\begin{lemma}\label{lemma: Gauss integral evaluation on the easy set}
Let $X$ be any good regular element of $\gF$, let $\GF_X$ be the
centralizer of $X$ in $\GF$ and let $x$ be an element of the
building for $\GF_X$ in $\GF$. Set $r=d_x(X)$. Let $\vphi$ be the
normalized characteristic function of the $\Mf_x$-orbit of
$\rho_{x,r}(X)$ in $\mf_{x,r}$. Then
\begin{equation}
    \forall  Y \in \gF_{x,-r},\qquad\i_{x,X}(Y) = \meas(\GF_x)\
    \hphi_{x,-r}(Y),
\end{equation}
where $\meas$ refers to the Haar measure appearing in the
definition of $\i_x$. (See Definition~\ref{definition: finite
Fourier transform} for $\hphi$ and Definition~\ref{definition:
inflation} for $\hphi_{x,-r}$).
\end{lemma}

\begin{proof}
By Definition~\ref{definition: Gauss integral}, $\i_{x,X}(Y) =
\int_{\GF_x} \kill{Y}{\Ad(g)X}\, dg$. Let $dk$ denote the quotient
measure on $\GF_x / \GF_{x,0^+}$, so
\begin{eqnarray*}
\i_{x,X}(Y)
    &=& \int\limits_{\GF_x} \kill{Y}{\Ad(g)X}\, dg\\
    &=&\int\limits_{\GF_x / \GF_{x,0^+}} \int\limits_{\GF_{x,0^+}}
\kill{\Ad(k)^{-1}Y}{\Ad(h)X}\, dh\, dk.
\end{eqnarray*}
Consider the terms $\Ad(h)X$ and $\Ad(k)^{-1}Y$. When $\epsilon>0$
is sufficiently small, Lemma~\ref{lemma: Adler 2} asserts the
existence of a diffeomorphism $e_{x,\epsilon}=e_{x,0^+} :
\gF_{x,0^+} \to \GF_{x,0^+}$.  Let $Z$ be the unique element of
$\gF_{x,0^+}$ such that $e_{x,0^+}(Z) = h$. Then $\Ad(h)X$ is an
element of the coset $X + \ad(Z) X + \gF_{x,r^+}$. Let $\t$ denote
the Lie algebra of $\GF_X$. From \cite[1.9.3]{Adler}, we have
\begin{equation}\label{eqn:perp decomp}
    \begin{array}{lll}
    \gF_{x,s} &= {\t}_{x,s}\oplus {\t}_{x,s}^\perp\\
    \gF_{x,s^+} &= {\t}_{x,s^+}\oplus {\t}_{x,s^+}^\perp
    \end{array},
\end{equation}
where ${\t}_{x,s}^\perp = \gF_{x,s} \cap {\t}^\perp$ and
${\t}_{x,s^+}^\perp = \gF_{x,s^+} \cap {\t}^\perp$. Together with
Lemma~\ref{lemma: Adler 1}, it follows that the coset $\ad(Z) X +
\gF_{x,r^+}$ is contained in $\gF_{x,r^+}$. Thus,
\begin{equation}
\Ad(h) X \in X + \gF_{x,r^+} \subset \gF_{x,r}.
    \label{eqn:hX}
\end{equation}
Notice also that
\begin{equation}
\Ad(k^{-1})Y \in \gF_{x,-r},
    \label{eqn:kX}
\end{equation}
since $k \in \GF_x$ and $Y\in \gF_{x,-r}$. Now argue as in the
proof of Lemma~\ref{proposition: inflation and FT}. From
Equations~\ref{eqn:hX} and \ref{eqn:kX}, it follows that
\begin{equation}
\kill{\Ad(k)^{-1}Y}{\Ad(h)X} =
\kill{\Ad(k)^{-1}Y}{X}.\label{eqn:khX}
\end{equation}

Now, return to $\i_{x,X}$ as above equipped with
Equation~\ref{eqn:khX} and we have
\begin{eqnarray*}
\i_{x,X}(Y)
    &=& \int\limits_{\GF_x / \GF_{x,0^+}} \int\limits_{\GF_{x,0^+}}
\kill{\Ad(k)^{-1}Y}{\Ad(h)X}\, dh\, dk\\
    &=& \int_{\GF_x / \GF_{x,0^+}}\kill{Y}{\Ad(k)X}\, dk
        \int_{\GF_{x,0^+}}\, dh\\
    &=& \meas(\GF_x) \vert \Mf_x\vert^{-1} \int_{\GF_x / \GF_{x,0^+}}
\kill{Y}{\Ad(k)X}\, dk\\
    &=& \meas(\GF_x) {\vert \Mf_x \vert}^{-1}\ \sum_{m \in \Mf_x}
\Lambda_{x,r}(\Ad(m)\X,\Y),
\end{eqnarray*}
where $\X = \rho_{x,r}(X)$. If $\orbit(\X)$ is the $\Mf_x$-orbit of $\X$ in
$\mf_{x,r}$ and if $1_{\orbit(\X)}$ denotes the characteristic function of
$\orbit(\X)$, then
\begin{eqnarray*}
&& \meas(\GF_x) {\vert \Mf_x \vert}^{-1}\ \sum_{m \in \Mf_x}
\Lambda_{x,r}(\Ad(m)\X,\Y)\\
 &=& \meas(\GF_x) \frac{\vert
Z_{\Mf_x}(\X) \vert}{\vert \Mf_x \vert}\
        \sum_{\Z \in \mf_{x,r}} \Lambda_{x,r}(\Z,\Y)\
        1_{\orbit(\X)}(\Z),
\end{eqnarray*}
where $Z_{\Mf_x}(\X)$ is the centralizer of $\X$ in $\Mf_x$.
Recognize the sum above as a Fourier transform on $\mf_{x,r}$, so
\begin{eqnarray*}
&&\meas(\GF_x) \frac{\vert Z_{\Mf_x}(\X) \vert}{\vert \Mf_x
\vert}\
        \sum_{\Z \in \mf_{x,r}} \Lambda_{x,r}(\Z,\Y)\
        1_{\orbit(\X)}(\Z)\\
         &=& \meas(\GF_x)\ {\vert \orbit(\X) \vert}^{-1}\
\widehat{1_{\orbit(\X)}}(\Y).
\end{eqnarray*}
Since ${\vert \orbit(\X) \vert}^{-1}\ 1_{\orbit(\X)}$ is the
function $\vphi$ defined in the statement of the lemma, then we
have
\begin{eqnarray*}
&& \meas(\GF_x)\ {\vert \orbit(\X) \vert}^{-1}\
\widehat{1_{\orbit(\X)}}(\Y)\\
    &=& \meas(\GF_x)\ \hphi(\rho_{x,-r}(Y))\\
    &=& \meas(\GF_x)\ \hphi_{x,-r}(Y).
\end{eqnarray*}
This proves Lemma~\ref{lemma: Gauss integral evaluation on the
easy set}.
\end{proof}

\begin{lemma}\label{lemma: vanishing}
Let $X$ be any good regular element of $\gF$, let $\GF_X$ be the
centralizer of $X$ in $\GF$ and let $x$ be an element of the
building for $\GF_X$ in $\GF$. Set $r=d_x(X)$. If $Y \in \gF_{-r}$
and $Y \not\in \gF_{x,-r}$ then $\i_{x,X}(Y) =0$.
\end{lemma}

\begin{proof}
Let $s = d_x(Y)$ and $t = -(r + s)/2$. Notice that $s < -r$ since
$Y \not\in \gF_{x,-r}$; thus $t > 0$. Let $\mathcal{G}_{x,t}$ be a
set of representatives for the $\GF_{x,t^+}$-cosets in $\GF_{x}$.
Then,
\begin{equation}\label{eqn:isum}
\i_{x,X}(Y) = \sum_{g_0 \in\, \mathcal{G}_{x,t} } \int_{\GF_{x,t^+}}
\kill{\Ad(g)X}{\Ad(g_0)Y}\, dg.
\end{equation}
Fix $g_0 \in \mathcal{G}_{x,t}$ and let $Y_0 = \Ad(g_0)Y$. Note
that $d_x(Y_0) = d_x(Y) =s$, as $Y_0$ is conjugate to $Y$ by an
element of $\GF_x$. Let $t^\prime$ be the smallest jump point for
$x$ greater than $t$. Since $t$ is positive and $t^\prime > t$ we
use Lemma~\ref{lemma: Adler 2} to define $e_{x,t^\prime} :
\gF_{x,t^\prime} \to \GF_{x,t^\prime}$, which we then use to
pull-back the measure $dg$ on $\GF_{x,t^+}$ to a measure $dZ$ on
$\gF_{x,t^+}$. Thus,
\begin{equation}\label{eqn:dgdZ}
\int_{\GF_{x,t^+}} \kill{\Ad(g)X}{Y_0}\, dg = \int_{\gF_{x,t^+}}
\kill{X}{\Ad(e(Z))Y_0}\, dZ.
\end{equation}
From Lemma~\ref{lemma: Adler 2}, we conclude that
\begin{equation}
\Ad(e(Z)) Y_0 \in Y_0 + \ad(Z) Y_0 + \gF_{x,2t^\prime + s}.
\end{equation}
Thus,
\begin{equation}\label{eqn:dZAd}
\Ad(e(Z)) Y_0 = Y_0 + \ad(Z) Y_0 + W,
\end{equation}
for some $W \in \gF_{x,2t^\prime + s}$. Since $2t+s + r =0$ (by the
definition of
$t$) and since $t^\prime > t$, it follows that  $2t^\prime + s + r > 0$.
Thus,
$\kill{X}{W}=1$ and
\begin{equation}\label{eqn:AddZ}
\kill{X}{\Ad(e(Z))Y_0} = \kill{X}{Y_0}\  \kill{X}{\ad(Z)Y_0}.
\end{equation}
Now combine Equations~\ref{eqn:dgdZ}, \ref{eqn:dZAd} and
\ref{eqn:AddZ} to conclude that
\begin{equation}
\int\limits_{\GF_{x,t^+}} \kill{\Ad(g)X}{Y_0}\, dg = \kill{X}{Y_0}
\int\limits_{\gF_{x,t^+}} \kill{X}{\ad(Z)Y_0}\, dZ.
\end{equation}

For a contradiction, suppose now that $\i_{x,X}(Y) \ne 0$. Fix
$Y_0 = \Ad(g_0)Y$ such that
\begin{equation}
\int_{\GF_{x,t^+}} \kill{\Ad(g)X}{Y_0}\, dg \ne 0.
\end{equation}
By Equation~\ref{eqn:isum}, this implies
\begin{equation}\label{eqn:ne0}
\int_{\gF_{x,t^+}} \kill{X}{\ad(Z)Y_0}\, dZ\ne 0.
\end{equation}
Now
\begin{eqnarray*}
 \kill{X}{\ad(Z)Y_0}
    &=& \kill{X}{[Z,Y_0]}\\
    &=& -\kill{X}{[Y_0,Z]}\\
    &=& -\kill{X}{\ad(Y_0)Z}\\
    &=& \kill{\ad(Y_0)X}{Z}\\
    &=& \kill{[Y_0,X]}{Z}\\
    &=& -\kill{[X,Y_0]}{Z}\\
    &=& -\kill{\ad(X)Y_0}{Z}.
\end{eqnarray*}
Thus, equation~\ref{eqn:ne0} becomes
\begin{equation}\label{eqn:ne1}
\int_{\gF_{x,t^+}} \kill{\ad(X)Y_0}{Z}\, dZ \ne 0.
\end{equation}

Now $Z \mapsto \kill{\ad(X)Y_0}{Z}$ defines an additive character
on $\gF_{x,t^+}$, so the integral above is nonzero if and only if
that character is trivial on $\gF_{x,t^+}$. We may therefore
assume that $\kill{\ad(X)Y_0}{Z} =1$ for all $Z\in \gF_{x,t^+}$,
and it follows that $\ad(X)Y_0 \in \gF_{x,-t}$. Equivalently, we
conclude that
\begin{equation}\label{eqn:depth geq -t}
d_x([X,Y_0]) + t \geq 0.
\end{equation}

To simplify notation, let $\t = \Lie(\GF_X)$ and let $\t^\perp$
denote subspace of $\gF$ perpendicular to $\t$ with respect to the
killing form. Write $Y_0 = Y_0^{\prime} + Y_0^\perp$ where
$Y_0^{\prime} \in {\t}_{x,s} \setminus {\t}_{x,s^+}$ and
$Y_0^\perp \in {\t}_{x,s}^\perp \setminus {\t}_{x,s^+}^\perp$,
using Equations \ref{eqn:perp decomp} again
(cf:\cite[1.9.3]{Adler}). Notice that $Y_0^\prime$ is semi-simple.
Now study $[X,Y_0] = [X,Y_0^\perp]$. From Lemma~\ref{lemma: Adler
1}, recall that $d_x([X,Y_0^\perp]) = d_x(X) + d_x(Y_0^\perp)$.
Thus, by Equation~\ref{eqn:depth geq -t}, $t \geq -r
-d_x(Y_0^\perp)$. This relation forces $Y_0^\perp$ deeper than
$Y_0$, since $d_x(Y_0^\perp) \geq -t -r = t + s > s$ (recall that
$t$ is strictly positive). Therefore, from $Y_0 = Y_0^{\prime} +
Y_0^\perp$, it follows that
\begin{equation}
Y_0 \in  Y_0^\prime + \gF_{x,s^+}.
\end{equation}

We have not yet used the fact that $Y \in \gF_{-r}$. By
\cite[3.1.2,~part~1]{Adler-DeBacker},
\begin{equation}
\gF_{-r} \subseteq \gF_\nilp+\gF_{x,-r}.
\end{equation}
So, we write
\begin{equation}\label{eqn:Y0}
Y_0 =  N + Y_0^{\prime\prime},
\end{equation}
where $Y_0^{\prime\prime} \in \gF_{x,-r}$ and $N \in \gF_{\nilp}$. Since $s
=
d_x(Y) = d_x(Y_0)$ and since $s < -r$, it follows that $Y_0^{\prime\prime}
\in
\gF_{x,s^+}$. Thus,
\begin{equation}\label{eqn:Y0N}
Y_0 \in N + \gF_{x,s^+}.
\end{equation}
Combine Equations~\ref{eqn:Y0} and \ref{eqn:Y0N} to see that
\begin{equation}
N \in Y_0^\prime + \gF_{x,s^+}.
\end{equation}
However, by Lemma~\ref{lemma: Adler 3}, the coset $Y_0^\prime +
\gF_{x,s^+}$ contains no nilpotent elements. This is the desired
contradiction and proves Lemma~\ref{lemma: vanishing}.
\end{proof}

\begin{proposition}\label{proposition: Gauss integral evaluation}
Let $X$ be any good regular element in $\gF$, let $\GF_X$ be the
centralizer of $X$ in $\GF$ and let $x$ be an element of the
building for $\GF_X$ in $\GF$. Set $r=d_x(X)$. Then
\begin{equation}
    \forall  Y\in \gF_{-r},\qquad\i_{x,X}(Y) = \meas(\GF_x)\
\hphi_{x,-r}(Y),
\end{equation}
where $\vphi$ is the normalized characteristic function of the $\Mf_x$-orbit
of
$\rho_{x,r}(X)$ in $\mf_{x,r}$.
\end{proposition}

\begin{proof}
If $Y$ is not an element of $\gF_{x,-r}$ then $\hphi_{x,-r}(Y)
=0$. With this observation, Proposition~\ref{proposition: Gauss
integral evaluation} follows immediately from Lemmas~\ref{lemma:
Gauss integral evaluation on the easy set} and \ref{lemma:
vanishing}.
\end{proof}


\subsection{Good, regular, elliptic orbital integrals}\label{section:
elliptic integrals}

This section uses Gauss integrals to produce an integral formula
for good, regular, elliptic orbital integrals.

\begin{definition}\label{definition: orbital integral}
For any $f\in \Sch{\gF}$ and $X\in \gF$, define
\begin{equation}
\orbint_\GF(X,f) = \int_{\GF_X\leftcoset\GF}f(\Ad(g^{-1}) X)\, dg.
\end{equation}
We will refer to $\orbint_\GF(X,f)$ as the {\em orbital integral
of $f$ at $X$}. Here, $dg$ refers to the quotient measure on
${\GF_X\leftcoset\GF}$.
\end{definition}

\begin{proposition}\label{proposition: elliptic}
Let $X$ be a tamely ramified good regular elliptic element of
$\gF$, let $\GF_X$ denote the centralizer of $X$ in $\GF$ and let
$x$ be an element of the building for $\GF_X$ in $\GF$. Set $r =
d_x(X)$. Let $\vphi\in\C(\mf_{x,r})$ be the normalized
characteristic function of the $\Mf_x$-orbit of $\rho_{x,r}(X)$ in
$\mf_{x,r}$, as in Proposition~\ref{proposition: Gauss integral
evaluation}. Then
\begin{equation}
\forall f\in \Sch{\gF}_r,
\qquad \orbint_\GF(X,f) =
\int_{\GF_X\leftcoset\GF} \int_{\gF} f(\Ad(g)^{-1}Y)\
\frac{\vphi_{x,r}(Y)}{\vol(\gF_{x,r^+})} \, dY\, dg.
\end{equation}
\end{proposition}

\begin{proof}
Let $\AG$ be the split component of the center of $\GF$ and let
$\TF$ denote $\GF_X$. Notice that the function $g\mapsto
f(\Ad(g)^{-1}X)$ factors from $\TF\leftcoset\GF$ to a function on
the left coset of $\GF\rightcoset\AG$ by $\TF\rightcoset\AG$.
Equip $\AG$ with a measure such that the compact set
$\TF\rightcoset\AG$ has measure $1$ with respect to the quotient
measure. Then,
\begin{equation}
\orbint_\GF(X,f) = \int_{\GF\rightcoset\AG} f(\Ad(g)^{-1}X)\, d^*g,
\end{equation}

Using Proposition~\ref{proposition: depth r Hecke} and
Definition~\ref{definition: depth r Hecke}, let $\check{f}$ be the
unique element of $\Sch{\gF}_{-r}$ such that $\widehat{\check{f}}
= f$. Then
\begin{equation}
\int_{\GF\rightcoset\AG} f(\Ad(g)^{-1}X)\, d^*g =
\int_{\GF\rightcoset\AG} \int_\gF
        \kill{Y}{\Ad(g)^{-1}X}\ \check{f}(Y)
                        \, dY\, d^*g.
\end{equation}
Pass from integration on $\GF\rightcoset\AG$ to integration on
$(\GF_x\leftcoset\GF\rightcoset\AG) \times \GF_x$; to simplify
notation, write $d_x^*$ for the quotient measure on
$\GF_x\leftcoset\GF\rightcoset\AG$. Then
\begin{eqnarray*}
&&
\int_{\GF\rightcoset\AG} \int_\gF
        \kill{Y}{\Ad(g)^{-1}X}\ \check{f}(Y)
                        \, dY\, d^*g\\
 &=& \int_{\GF_x\leftcoset\GF\rightcoset\AG} \int_{\GF_x}
\int_{\gF}
        \kill{\Ad(kg)^{-1}X}{Y}\ \check{f}(Y)
                \, dY\, dk\, d_x^*g.
\end{eqnarray*}
A change of variables gives
\begin{eqnarray*}
&&    \int_{\GF_x\leftcoset\GF\rightcoset\AG} \int_{\GF_x}
\int_{\gF}
        \kill{\Ad(kg)^{-1}X}{Y}\ \check{f}(Y)
                \, dY\, dk\, d_x^*g\\
    &=& \int_{\GF_x\leftcoset\GF\rightcoset\AG}
    \int_{\GF_x} \int_{\gF} \kill{\Ad(k)^{-1}X}{Y}\ \,
    \check{f}(\Ad(g)^{-1}Y) dY\, dk\, d_x^*g.
\end{eqnarray*}
Since $\check{f}$ is supported by a compact set, the last two
integrals may be exchanged, giving
\begin{eqnarray*}
&& \int_{\GF_x\leftcoset\GF\rightcoset\AG}
    \int_{\GF_x} \int_{\gF} \kill{\Ad(k)^{-1}X}{Y}\ \,
    \check{f}(\Ad(g)^{-1}Y) dY\, dk\, d_x^*g\\
     &=& \int\limits_{\GF_x\leftcoset\GF\rightcoset\AG}
\int\limits_{\gF} \int\limits_{\GF_x} \kill{\Ad(k)X}{Y}\, dk\
\check{f}(\Ad(g)^{-1}Y)\, dY\, d_x^*g.
\end{eqnarray*}
The integral over $\GF_x$ is $\i_x(X,Y)$, so we have shown
\begin{equation}
\orbint_\GF(X,f) =
        \int_{\GF\rightcoset\AG}\int_{\gF}
         \i_x(X,Y)\ \meas(\GF_x)^{-1}\ \check{f}(\Ad(g)^{-1}Y)
                \, dY\, d^*g.
\end{equation}

Now we are in a position to make use of the results from
Section~\ref{section: Fourier transform}; namely, since $x$ is a
point in the building for $\TF$ in $\GF$, then by
Proposition~\ref{proposition: Gauss integral evaluation},
$\i_x(X,Y) = \meas(\GF_x)\ \hphi_{x,-r}(Y)$, for any $Y$ in
$\gF_{-r}$. Since the support of $\check{f}$ is contained in
$\gF_{-r}$, it follows that
\begin{eqnarray*}
&&\int_{\GF\rightcoset\AG}\int_{\gF}
         \i_x(X,Y)\ \meas(\GF_x)^{-1}\ \check{f}(\Ad(g)^{-1}Y)
                \, dY\, d^*g\\
    &=& \int_{\GF\rightcoset\AG} \int_{\gF}
        \hphi_{x,-r}(Y)\ \check{f}(\Ad(g)^{-1}Y)\, dY\, d^*g.
\end{eqnarray*}
Using Proposition~\ref{proposition: inflation and FT}, we have
\begin{eqnarray*}
&& \int_{\GF\rightcoset\AG} \int_{\gF}
        \hphi_{x,-r}(Y)\ \check{f}(\Ad(g)^{-1}Y)\, dY\, d^*g\\
    &=&
        \int_{\GF\rightcoset\AG} \int_{\gF}
\frac{\widehat{\vphi_{x,r}}(Y)}{\vol(\gF_{x,r^+})}\check{f}(\Ad(g)^{-1}Y)\,
dY\, d^*g\\
    &=&    \int_{\GF\rightcoset\AG} \int_{\gF}
\frac{\vphi_{x,r}(Y)}{\vol(\gF_{x,r^+})}\widehat{\check{f}}(\Ad(g)^{-1}Y)\,
dY\, d^*g.
\end{eqnarray*}
Since $f = \widehat{\check{f}}$, we have shown that
\begin{equation}
\orbint_\GF(X,f) =
        \int_{\GF\rightcoset\AG} \int_{\gF}
        f(\Ad(g)^{-1}Y) \frac{\vphi_{x,r}(Y)}{\vol(\gF_{x,r^+})}\, dY\,
d^*g.
\end{equation}
Recasting this result, observe that $\TF\rightcoset\AG$ is
represented by elements of $\GF_x$, so
\begin{equation}
\orbint_\GF(X,f) =
        \int_{\TF\leftcoset\GF} \int_{\gF}
        f(\Ad(g)^{-1}Y)\ \frac{\vphi_{x,r}(Y)}{\vol(\gF_{x,r^+})}\, dY\, dg,
\end{equation}
where $dg$ refers to the quotient measure on $\TF\leftcoset\GF$ suitable normalized.
This proves Proposition~\ref{proposition: elliptic}.
\end{proof}

\begin{example}
Suppose $\GF = SL(2,\F)$ and $f = 1_{\gF(\A)}$; let
$\orbint_\gF(X)$ denote $\orbint_\GF(X,f)$.  We now illustrate
Proposition~\ref{proposition: elliptic} and foreshadow ideas
presented in Section~\ref{section: conjugacy}.  We fix a standard
fundamental affine chamber in the building for $\GF$ with
polyvertices $0$ and $1$; the facets of this chamber will be
denoted $(0)$, $(1)$ and $(01)$ with $\GF_{(0)} = SL(2,\A)$.
First, suppose $X$ is an element of the Cartan subalgebra
\[
\left\{
\begin{pmatrix}
  0 & a \\
  \varepsilon a & 0
\end{pmatrix}
\ \vert\
a\in F
\right\}.
\]
where $\varepsilon$ is a fixed quadratic residue in the group
$\A^*$ of units of $\A$. Then the image of $\B(\GF_X)
\hookrightarrow \B(\GF)$ is ${(0)}$; let $x = (0)$. The depth of
$X$ is $0$ if and only if $a\in \A^*$; suppose that is the case
and set $r=0$. Now $\bar{\gF}_{x,r} = sl(2,\f)$ and this is
equipped with the adjoint action of $\bar{\GF}_x = SL(2,\f)$. Then
the $\bar{\GF}_x$-orbit $\orbit_{\bar{GF}_x}(\rho_{x,r}(X))$ of
$\rho_{x,r}(X)\in \bar{\gF}_{x,r}$ is the set of $\f$-rational
points on the smooth variety
\[
\left\{
\begin{pmatrix}
  z & x \\
  y & -z
\end{pmatrix}
\ \vert\
z^2+xy = \varepsilon \bar{a}^2
\right\},
\]
where $\bar{a}$ is the image of $a$ under $\A \to \f$. (Note that
$\bar{a} \ne 0$.) In particular, as $X$ varies in the set of depth
$0$ elements of this Cartan subalgebra, the value of
$\orbint_\gF(X)$ is determined by $\bar{a}$. Next, suppose $Y$ is
an element of the Cartan subalgebra
\[
\left\{
\begin{pmatrix}
  0 & b \\
  \varpi b & 0
\end{pmatrix}
\ \vert\
b\in F
\right\}.
\]
where $\varpi$ is a fixed uniformizer of $F$. Let $y\in B(\GF)$ be
the barycenter of the maximal facet $(01)$ of $\B(\GF)$. The depth
of $Y$ is $\frac{1}{2}$ if and only if $b \in \A^*$; suppose that is
the case and set $s= \frac{1}{2}$. Then $\bar{\gF}_{y,s} =
\mathbb{A}^2(\f)$ is equipped with the action of $\bar{\GF}_y =
GL(1,\F)$ given by $t\cdot(u,v) = (t^2u,t^{-2}v)$. Although we do
not pursue this point of view here, we remark that the
characteristic function of $\bar{\GF}_y$-orbit
$\orbit_{\bar{GF}_y}(\rho_{y,s}(Y))$ of $\rho_{y,s}(Y)\in
\bar{\gF}_{y,s}$ is the characteristic function of the
Frobenius-stable Kummer local system (an $\ell$-adic sheaf)
corresponding to the trivial character of the component group of
\[
\left\{
\begin{pmatrix}
  u, & v
\end{pmatrix}
\ \vert\
uv = \bar{b}^2
\right\},
\]
where $\bar{b}$ is the image of $b$ under $\A \to \f$. (Note that
$\bar{b} \ne 0$.) In particular, as $Y$ varies in the set of depth
$\frac{1}{2}$ elements of this Cartan subalgebra, the value of
$\orbint_\gF(Y)$ is determined by $\bar{b}$. The phenomena
illustrated by these examples will be generalized considerably in
Section~\ref{sec:parameter-space} and Theorem~\ref{thm:kappa}.
\end{example}


\subsection{Descent}\label{section: descent}

In this section we apply standard parabolic descent arguments to
extend Proposition~\ref{proposition: elliptic} to good, regular
orbital integrals; the result is Theorem~\ref{theorem: main},
which is the main result of Section~\ref{section: formula}, which
will allow use to compare orbital integrals without evaluating
them.

Here, $K$ is a fixed maximal special parahoric subgroup of $\GF$.

\begin{definition}\label{definition: Jacquet}
Let $\LF$ be a Levi subgroup of $\GF$ and let $\lF$ denote the Lie
algebra for $\LF$. Let $\PF$ be any parabolic subgroup of $\GF$
with Levi component $\LF$. Let $\UF$ be the unipotent radical of
$\PF$ and let $\uF$ denote the Lie algebra for $\UF$. For any
$f\in \Sch{\gF}$ define $f_\PF \in \Sch{\lF}$ by
\begin{equation}
    \forall Y \in \lF,\qquad f_\PF(Y) = \int_{\uF} \int_K
    f(\Ad(k)^{-1}(Y+Z))\, dk\, dZ.
\end{equation}
\end{definition}

\begin{remark}
The order of integration above is unimportant, as all relevant
integrands have compact support.
\end{remark}

\begin{theorem}\label{theorem: main}
Let $X$ be a good regular element of $\gF$ of depth $r$. Let $\LF$
be a Levi subgroup of $\GF$ containing $\GF_X$ as an elliptic
Cartan subgroup. Let $x$ denote the image of the building for
$\GF_X$ in $\LF$. Let $\vphi \in \C(\bar{\lF}_{x,r})$ be the
normalized characteristic function of the $\bar{\LF}_x$-orbit of
the image of $X$ under $\lF_{x,r} \to \lF_{x,r}/\lF_{x,r^+}$.
Then, for all  $f\in \Sch{\gF}_r$,
\begin{equation}
D^{\gF,\lF}(X)\
\orbint_\GF(X,f) = \int_{\GF_X\leftcoset\LF} \int_{\lF}
f_\PF(\Ad(g)^{-1}Y)\ \frac{\vphi_{x,r}(Y)}{\vol(\lF_{x,r^+})} \,
dY\, dg,
\end{equation}
where $D^{\gF,\lF}(X) =
\abs{\det(\ad(X)\vert_{\gF/\lF})}_\F^{1/2}$.
\end{theorem}

\begin{proof}
It is important to notice that $\vphi_{x,r}$ denotes a locally
constant function on $\lF$ supported by $\lF_{x,r}$, since $\vphi$
is an element of $\C(\bar{\lF}_{x,r})$.

From \cite{Harish-Chandra} we have
\begin{equation}\label{equation: orbital integral descent}
    \forall f\in\Sch{\gF},
\qquad
D^{\gF,\lF}(X)\ \orbint_\GF(X,f)
    = \orbint_\LF(X,f_\PF),
\end{equation}
where $X\in \lF$ is elliptic. From the definition of the depth
function (cf: Section \ref{section: preliminaries}) it follows
that the depth $r$ of $X\in\gF$ equals the depth of $X\in\lF$
relative to $x\in \B(\LF)$; thus, $r = d_x(X)$. We claim that
$f_\PF \in \Sch{\lF}_r$. To see this, we must show that
$\Fourier_\lF f_\PF\in \Sch{\lF}$ is supported by $\lF_{-r}$. From
\cite{Harish-Chandra} we recall that the Fourier transform
commutes with the operator $f \mapsto f_\PF$; more precisely, for
any $f\in C^\infty_c(\gF)$,
\begin{equation}
\Fourier_\lF(f_P) = (\Fourier_\gF f)_P.
\end{equation}
Recall that $f\in \Sch{\gF}_r$ if and only if $\Fourier_\gF f \in
\Sch{\gF_{-r}}$, by Definition~\ref{definition: depth r Hecke}.
From Definition~\ref{definition: Jacquet} we see that
$(\Fourier_\gF f)_P \in \Sch{\gF_{-r}\cap \lF}$. By
\cite[3.5.3]{Adler-DeBacker}, $\lF_{-r} = \gF_{-r} \cap \lF$,
which shows that $f_\PF \in \Sch{\lF}_r$. Now, by
Proposition~\ref{proposition: elliptic},
\begin{equation}\label{equation: orbital integral descended}
\orbint_\LF(X_\lF,f_\PF) =  \int_{\GF_X\leftcoset\LF} \int_{\lF}
f_\PF(\Ad(g)^{-1}Y)\ \frac{\vphi_{x,r}(Y)}{\vol(\lF_{x,r^+})} \,
dY\, dg.
\end{equation}
Combining Equations~\ref{equation: orbital integral descent} and
\ref{equation: orbital integral descended} proves
Theorem~\ref{theorem: main}.
\end{proof}


\subsection{Local constancy of good regular orbital
integrals}\label{local constancy}

In this section we use Theorem~\ref{theorem: main} to describe the
local constancy of $X \mapsto \orbint_\GF(X,\cdot)$, where
$\orbint_\GF(X,\cdot)$ is restricted to spaces of functions which
are relevant to the remainder of this paper.

\begin{definition}\label{4.2}
Let $\tF = \Lie(\TF)$ be a tamely ramified Cartan subalgebra of
$\gF$. Choose a point $y_\TF$ in the building for $\TF$ in $\GF$.
For each real number $s$ and for each $\Mf_{y_\TF}$-orbit $\orbit$
in $\mf_{y_{\TF},s}$, define
\begin{equation}
\tF_\orbit = \tF \cap d_{y_\TF}^{-1}(s) \cap
\rho_{y_\TF,s}^{-1}(\orbit).
\end{equation}
\end{definition}

\begin{lemma}\label{4.3}
Let $\TF$, $\tF$ and $y_\TF\in \B(\TF)$ be as in
Definition~\ref{4.2}. Then
\[
\{ \tF_\orbit \tq \orbit \subset \mf_{y_\TF,s}\}
\]
defines a partition of $\tF$, where $s$ ranges over $\R$ and
$\orbit$ ranges over all $\Mf_{y_\TF}$-orbits in $\mf_{y_\TF,s}$.
By restriction, this defines a partition of the set $\gF^\ellip$ of
tamely ramified regular elliptic elements in $\gF$.
\end{lemma}

\begin{proof}
$\tF$ is filtered in $s$ by $\tF_s = \tF \cap \gF_{y_\TF,s}$, so a
partition of $\tF$ is given by the sets $\tF \cap
d_{y_\TF}^{-1}(s)$. This partition is refined according to the
partition of $\mf_{y_\TF,s}$ into $\Mf_{y_\TF}$-orbits. Since this
is the partition above, Lemma~\ref{4.3} is proved.
\end{proof}

\begin{proposition}\label{4.4}
Let $\gF^{g,e}_0$ denote the set of good regular elliptic elements
in $\gF_0$. The function
\begin{eqnarray*}
\gF^{g,e}_0 &\to& \Sch{\gF}_0^\ast\\
X &\mapsto&\orbint_\GF(X,\, \cdot\, )
\end{eqnarray*}
is constant on the partition of $\gF^{g,e}_0$ defined by
restricting the partition of $\gF^\ellip$ given in
Lemma~\ref{4.3}.
\end{proposition}

\begin{proof}
Let $\orbit$ be an $\Mf_x$-orbit in $\mf_{x,r}$. For {\it any} $X$
in $\t_\orbit$, the function $\vphi \in \C(\mf_{x,r})$ appearing
in Proposition~\ref{proposition: elliptic} is the normalized
characteristic function of $\orbit$. This proves
Proposition~\ref{4.4}.
\end{proof}

\begin{corollary}\label{cor:unit element}
Let $f$ be the characteristic function of $\gF(\A)$ in
$\Sch{\gF}$. Let $X$ be a good regular element of $\gF$ with
non-negative depth $r$. Let $\LF$ and $x$ be as in
Theorem~\ref{theorem: main}. Write $\tF$ for the Cartan subalgebra
in $\gF$ for $\GF_X$, so $X\in \tF_{x,r}$. If $X^\prime \in
\tF_{x,r}$ is regular and $D^{\gF,\lF}(X) =
D^{\gF,\lF}(X^\prime)$, then
\begin{equation}
\rho_{x,r}(X) = \rho_{x,r}(X^\prime) \qquad\implies\qquad
\orbint_\GF(X,f) =\orbint_\GF(X^\prime,f).
\end{equation}
\end{corollary}

\begin{proof}
This follows from Theorem~\ref{theorem: main} and the following
facts: $f$ is an element of $\Sch{\gF}_r$ by Remark~\ref{remark:
unit}; the depth of $X$ in $\lF$ equals the depth of $X$ in $\gF$,
where $\lF$ is the Lie algebra for $\LF$; $X^\prime$ is good; and
for any parabolic $\PF$ with Levi component $\LF$, $f_\PF$ is the
characteristic function of $\lF(\A)$ which is an element of
$\Sch{\lF}_0$ and therefore of $\Sch{\lF}_r$.
\end{proof}


\section{Statement of Results}\label{sec:statement}

Proposition~\ref{4.4} and Corollary~\ref{cor:unit element} give
explicit results about the local constancy of orbital integrals.
The rest of this paper draws some implications from this formula
in the special case that $\lie{g}$ is a classical Lie algebra and
the function $f$ is the characteristic function of the
integral-valued points of $\lie{g}$.  We assume for the rest of
the paper that $f$ is that function.

In this special case, arithmetic motivic integration presents the
orbital integrals as the number of points on varieties over finite
fields. This presentation is independent of the underlying local
field in a sense that we will make precise below.

From the field-independent description, we deduce that the
fundamental lemma holds for a restricted set of elements for local
fields in zero characteristic, if the corresponding statement is
known in positive characteristic.

\subsection{Notation}

Recall that $F$ is a $p$-adic field with ring of integers $\A$,
prime ideal $\pFA$, and residue field $\f$, with $q=q_F$. Let
$\bar F$ and $\barf$ be the algebraic closures of $F$ and $\f$.
Let $\varpi = \varpi_F$ be a uniformizer in $F$.  We normalize the
absolute value so that $|\varpi_F|=q^{-1}$.  We extend the
normalized absolute value to an absolute value on $\bar F$. Let
$\op{res}:\A\to \f$ be the residue map.  We let $\op{ord}:\bar
F\to \ring{Q}$ be the valuation, normalized so that
   \begin{equation}
   |x| = q^{-\op{ord}\, x}.
   \end{equation}
Let
   \begin{equation}
   F^{int} = \{x \in \bar F \tq \op{ord}(x) \in \ring{Z}\}.
   \end{equation}
We let $\op{ac}:F^{int} \to \barf^\times$ be the angular component
function given by
   \begin{equation}
   \op{ac}(0) = 0,\quad \op{ac}(x) =
   \op{res}(x/\varpi^{\op{ord}\,x}).
   \end{equation}
It depends on a choice of uniformizer $\varpi$. For $i\in
\ring{Z}$, we let $\op{res}_i: F\to \f$ be the map
   \begin{equation}\op{res}_i(x) = \begin{cases} \ac{x} &
   \text{if}\op{ord}(x)=i\\
   0 & \text{otherwise}.
   \end{cases}
   \end{equation}

\begin{definition}
We say that a statement $\psi^F$ about local fields $F$ holds when
the residual characteristic is {\it sufficiently large\/} when
there is a natural number $M$ such that the statement holds
whenever $M$ is relatively prime to the characteristic of the
residue field of $F$. That is,
   \begin{equation}
   \exists M\ \forall F.\quad
   (q_F,M)=1 ~~\Rightarrow~~ \psi^F.
   \end{equation}
\end{definition}

Recall from Section~\ref{section: preliminaries} that we restrict
$F$ so that its residual characteristic is sufficiently large (for
the various statements that we make). This assumption is mentioned
in many of the lemmas and theorems, but even when it is not
explicitly mentioned, {\it the assumption remains in effect}.  The
natural number $M$ will depend on the Lie algebra $\lie{g}$ under
consideration and a parameter $r\in \ring{Q}$. For each $\lie{g}$
and $r$, the constant $M$ will be effectively computable.

\subsection{Lie algebras considered}

\begin{definition}{\label{def:pairs}}
Let $(\lie{g},\lie{h})$ be one of the following pairs of Lie
algebras.
    \begin{equation}
    \begin{array}{lll}
    \lie{so}(2c+1),\quad \lie{so}(2a+1) \oplus \lie{so}(2b+1), \text{ with }
a+b=c,\\
    \lie{sp}(2c),\quad \lie{sp} (2a)\oplus \lie{so}(2b), \text{ with } a+b =
c,\quad (b\ne1)\\
    \lie{so}(2c),\quad \lie{so}(2a)\oplus \lie{so}(2b),\text{ with } a + b =
c,\quad (a\ne1,b \ne1, c\ne1)\\
    \end{array}
    \end{equation}
\end{definition}

We refer to these three cases as the odd orthogonal, symplectic,
even orthogonal respectively.  In each case, the Lie algebra
$\lie{h}$ is a sum of two factors $\lie{h}_1\oplus \lie{h}_2$. We
write $(X,Y)$ for an element in $\lie{g}\oplus\lie{h}$, with $Y=
(Y_1,Y_2)\in\lie{h}_1\oplus\lie{h}_2$.  Below, we will fix a
concrete representation (the standard representation) of these
algebras.

These pairs are considered in the paper \cite{VTF}.  In that
paper, an additional family $\lie{u}(c)$, the Lie algebra of the
unitary group, is considered.  We make a few comments about the
unitary case in Remark~\ref{remark:unitary}.

\begin{remark}
The origin of this list of Lie algebras is the following. Let $G$
be a classical split adjoint group over $F$ and let $H$ be an
elliptic endoscopic group of $G$. Then the Lie algebras
$\lie{g},\lie{h}$ listed above are the Lie algebras of $G$ and
$H$. The list is not exhaustive.  In particular, it only includes
split cases.  We refer to $\lie{h}$ as an endoscopic algebra.
\end{remark}

For a given $\lie{g}$, if the residual characteristic of $F$ is
sufficiently large,  every Cartan subalgebra of
$\lie{g}\oplus\lie{h}$ splits over a tamely ramified extension of
$F$.  We confine our attention to local fields $F$ for which this
is the case.

Each of the lie algebras under consideration comes with a natural
representation.  We identify elements of $\gF$ with the matrices
that represent them.  We take the eigenvalues $\lambda(X)$ of a
semi-simple element with respect to this representation.

\begin{definition}\label{def:vg}
Let $\lie{g}$ be one of the semi-simple Lie algebras introduced in
Definition~{\ref{def:pairs}}. We say that an element $X$ is {\it
restricted } (of \emph{slope} $r\in\ring{Q}$) if it satisfies the
following conditions.
   \begin{enumerate}
   \item $X$ is regular semi-simple.
   \item $X$ is contained in a tamely ramified Cartan subalgebra $\tF$.
   \item $|\alpha(X)|=q^{-r}$ for each (absolute) root of $\gF$ relative to
$\tF$.
   \item $|\lambda(X)| = q^{-r}$, for each nonzero eigenvalue
   $\lambda$.
   \item The multiplicity of the eigenvalue $\lambda=0$ is at most $1$.
   \end{enumerate}
Write $\liegood{g}(r)$ for the set of restricted  elements of
slope $r$ in the lie algebra $\lie{g}$.  When it becomes necessary
to indicate the coefficient ring $A$ of the matrices $X$, we write
$\liegood{g}(r,A)\subset \lie{g}(A)$.
\end{definition}

The first three conditions in the definition of restricted imply
that every restricted element is good.  On the set of regular
elements in a Cartan subalgebra, the depth is equal to the slope
up to a nonzero multiplicative factor that depends only on the
Cartan subalgebra.  Although it is possible to give a formula for
this multiplicative factor, our proof does not rely on the value
of this scalar.  In the interest of simplicity, we do not give a
formula.

In the symplectic and odd orthogonal algebras, the final two
conditions follow from the first three, at least for fields of
sufficiently large residual characteristic (which we assume).
Finally, in the even orthogonal lie algebras, it is possible to
satisfy the first three conditions without the last two
conditions, because of a pair of eigenvalues $\pm\lambda(X)$ that
have smaller absolute value than the rest. (For instance, a good
element may have a pair of eigenvalues equal to zero.)

\begin{definition}\label{def:equiv}
We give an equivalence relation on the restricted elements of
slope $r$.  If the lie algebra is symplectic or odd orthogonal, we
say that two restricted  elements $X$ and $X'$ of slope $r$ are
{\it equivalent\/} if the eigenvalues $\lambda_i(X)$ of $X$ and
$\lambda_i(X')$ of $X'$ can be indexed so that
   \begin{equation}
   |\lambda_i(X')-\lambda_i(X)|< q^{-r},
   \label{eqn:alpha-equiv}
   \end{equation}
for all $i$.  For the even orthogonal lie algebra, let $J$ be the
symmetric matrix that defines the algebra:
    $$
    \lie{so}(2c) = \{X \tq {}^tX J + J X = 0\}.
    $$
We say that $X$ and $X'$ are equivalent if the
Inequality~\ref{eqn:alpha-equiv} holds and if the additional
condition
   \begin{equation}
   |\op{pfaff}(JX) - \op{pfaff}(J X')| < q^{-c r},
   \end{equation}
where $\op{pfaff}$ is the pfaffian\footnote{Pfaffians are
discussed further in Section~\ref{sec:pfaff}.} of a skew-symmetric
matrix. If $X$ is a restricted  element of slope $r$, let $[X]_r$
be its equivalence class.
\end{definition}

\begin{theorem} \label{thm:equiv-class}
Let $\lie{g}$ be one of the lie algebras given in
Definition~\ref{def:pairs}.  There exists $M>0$ and an affine
variety $S_{\lie{g},r}$ over $\ring{Z}[\frac{1}{M}]$ that
classifies the equivalence classes of restricted  elements of
slope $r$. The variety $S_{\lie{g},r}$ depends on $\lie{g}$ and
$r$, but is independent of the local field in the following sense.
For all local fields $F$ whose residual characteristic is prime to
$M$, we have a natural bijection between
  \begin{equation}
  \{[X]_r \tq X\in \liegood{g}(r) \}
  \end{equation}
and $S_{\lie{g},r}(\f)$, where $\f$ is the residue field of $F$.
\end{theorem}

The varieties $S_{\lie{g},r}$ are described for each $\lie{g}$ and
$r$ in Section~\ref{sec:parameter-space}.

\begin{proof} This will be proved later as Theorem~\ref{thm:bij}.
\end{proof}

\subsection{Orbital Integrals}

Langlands and Shelstad attach a $\kappa$-orbital integral to
semi-simple elements in the endoscopic algebra whose {\it image\/}
in $\lie{g}$ is regular. The $\kappa$-orbital integral of $f$ (the
characteristic function of $\lie{g}(\A)$) is defined in terms of a
transfer factor defined in \cite{Xf}.  We write the
$\kappa$-orbital integral -- including the transfer factor -- on
the Lie algebra, rather than the group. We use the canonical
normalization of transfer factors from \cite{HS}. When
$Y\in\lie{h}(r,\A)$, write
   \begin{equation}\orb{g}{h}(Y)\end{equation}
for the $\kappa$-orbital integral on $\lie{g}$ attached to $Y$
over the characteristic function of $\lie{g}(\A)$. It is a sum of
orbital integrals in $\lie{g}$, weighted by the Langlands-Shelstad
transfer factor. We write the stable orbital integral of the
characteristic function of $\lie{h}(\A)$ on $\lie{h}$ as
$\orb{h}{h}(Y)$.  If $Y$ is a regular semi-simple element of
$\lie{h}$, then Langlands and Shelstad have defined a
corresponding element $X\in\lie{g}$, which they call the {\it
image\/} of $Y$. The image of $Y$, which will be made explicit in
Section~\ref{section:kappa}, is well-defined up to stable
conjugacy.

\begin{conjecture}[The fundamental lemma]
For each of the pairs $\lie{g},\lie{h}$ in
Definition~\ref{def:pairs}, and every regular semi-simple element
$Y$ in $\lie{h}$, such that the image of $Y$ in $\lie{g}$ is
restricted of slope $r$, there is an equality of orbital integrals
   \begin{equation}
   \orb{g}{h}(Y) = \orb{h}{h}(Y).
   \label{eqn:fl}
  \end{equation}
\end{conjecture}
\smallskip

The endoscopic lie algebras that we study are given as a sum
  \begin{equation}\lie{h} = \lie{h}_1\oplus \lie{h}_2.\end{equation}
Each restricted  element $Y$ in $\lie{h}$ is an ordered pair
$(Y_1,Y_2)$ of two restricted  elements of the same slope. By
Theorem~\ref{thm:equiv-class}, the equivalence classes of
restricted elements of slope $r$ in $\lie{h}$ are parameterized by
$S_{\lie{h}_1,r}\times S_{\lie{h}_2,r}.$  We write this product of
parameter spaces as $S_{\lie{h},r}$.

If $Y$ and $Y'$ are equivalent and are restricted of slope $r$ in
$\lie{h}$, then the image of $Y$ in $\lie{g}$ is restricted of
slope $r$ if and only if the image of $Y'$ has the same property.
There is a subvariety $S_{\lie{g},\lie{h},r}$ of $S_{\lie{h},r}$
that parameterizes equivalence classes of elements $Y$ whose image
in $\lie{g}$ is restricted of slope $r$.

\begin{theorem}  Assume that the residual
characteristic is sufficiently large. Let $Y\in\lie{h}$. Assume
that $Y$ is a restricted  element of slope $r\ge0$.  The
$\kappa$-orbital integral
  \begin{equation}\orb{g}{h}(Y)\end{equation}
depends only on the equivalence class of $Y$.
\end{theorem}

\begin{proof}  This will be proved later as
Theorem~\ref{thm:kappa}.
\end{proof}

Consequently, we may speak of the $\kappa$-orbital integral
   \begin{equation}
   \orb{g}{h}(y)
   \end{equation}
of a parameter $y\in S_{\lie{g},\lie{h},r}(\f)$ when the
characteristic of $\f$ is sufficiently large.

If $U$ is a variety over any base variety $S$ and $x$ is a closed
point of the base $S$ with residue field $k(x)=\f$, then we let
$|U_x(\f)|$ be the number of points of the fiber of $U$ over $x$.
The following is the main result of the paper.

\begin{theorem}\label{thm:variety}
For every $(\lie{g},\lie{h})$ in Definition~\ref{def:pairs} and
$r\in\ring{Q}$ with $r\ge0$, there are a natural number $M$, a
finite indexing set $I$, varieties $U_i$ over
$S_{\lie{g},\lie{h},r}$ indexed by $i\in I$, constants
$b_i\in\ring{Q}$ indexed by $i\in I$, and a polynomial $p(x)$ of
the form
   \begin{equation}
   p(x) = x^k \prod_{i=1}^{k'} (x^{k_i} - 1),
   \end{equation}
with the following property: For all finite fields of order
relatively prime to $M$, we have
  \begin{equation}
  \forall y\in S_{\lie{g},\lie{h},r}(\f),\quad \orb{g}{h}(y)
    = \frac{1}{p(q)}\sum_{i\in I} b_i |U_{i,y}(\f)|.
  \end{equation}
The constant $M$, the indexing set $I$, the varieties $U_i$, and
the constants $b_i$, and the polynomial $p(x)$  are effectively
computable.
\end{theorem}

\begin{proof}  This will be proved later.  It is a consequence of
the fact that the integrals in question can be described as
volumes of a family of locally constant definable sets
(Lemma~\ref{lemma:psi-Pas}, Lemma~\ref{lemma:phi-stable},
Lemma~\ref{lemma:fld}) and that families of volumes of locally
constant definable sets have a representation of this form
(Theorem~\ref{thm:existU}).
\end{proof}

The theorem asserts that the $\kappa$-orbital integrals of
restricted elements of slope $r$ count points on varieties $U_i$
over finite fields.   We emphasize that the varieties $U_i$ and
constants $b_i$ are independent of the local field $F$ and the
residue field $\f$. These are universal varieties that calculate
the $\kappa$-orbital integrals for all local fields with
sufficiently large residual characteristic. Whenever it becomes
necessary to indicate the dependence of the data $I$, $U_i$,
$b_i$, and $p$ on the underlying parameters $(\lie{g},\lie{h})$
and $r$, we write
   \begin{equation}I=I(\lie{g},\lie{h},r),\
U_i=U(\lie{g},\lie{h},r)_i,\end{equation}
and so forth.

\begin{corollary}  For local fields of sufficiently large residual
characteristic, the fundamental lemma holds for restricted
elements of slope $r$ in $(\lie{g},\lie{h})$ iff
Equation~\ref{eqn:gfl} holds.
  \begin{equation}
  \begin{array}{lll}
  &\forall y\in S_{\lie{g},\lie{h},r}(\f).\quad\\
  &\quad\frac{1}{p(\lie{g},\lie{h},r)(q)}\sum_{i\in
  I(\lie{g},\lie{h},r)}
   b(\lie{g},\lie{h},r)_i
   |U(\lie{g},\lie{h},r)_{i,y}(\f)| =\\
  &\quad\frac{1}{p(\lie{h},\lie{h},r)(q)} \sum_{i\in I(\lie{h},\lie{h},r)}
   b(\lie{h},\lie{h},r)_i
   |U(\lie{h},\lie{h},r)_{i,y}(\f)|.
   \end{array}
   \label{eqn:gfl}
  \end{equation}
\end{corollary}

\begin{remark}   The varieties on the left are geometrizations of
the  $\kappa$-orbital integrals.  Those on the right are the
geometrizations of the stable orbital integrals. The stable
orbital integrals on $\lie{h}$ of the characteristic function of
the unit lattice $\lie{g}(\A)$ are products
  \begin{equation}
  \orb{h}{h}(Y) = \orb{h_1}{h_1}(Y_1)\times
  \orb{h_2}{h_2}(Y_2).
  \end{equation}
Thus, we may apply the results of Theorem~\ref{thm:variety} twice,
once for $\lie{h}_1$ and once for $\lie{h}_2$ to get a
representation of the stable orbital integral on $\lie{h}$ as the
number of $\f$-points on a formal linear combination of varieties.
It is this combination that appears on the right-hand side of the
corollary.)
\end{remark}

The representation of orbital integrals in
Theorem~\ref{thm:variety} is independent of the field $F$.  This
observation leads to the following corollary.

\begin{corollary} \label{cor:F}
If two local fields $F$, $F'$ (of sufficiently large residual
characteristic) have the same residue field $\f$, if $Y$ is
restricted  of slope $r$ in $\liegood{h}(r,F)$ and $Y'$ is
restricted of slope $r$ in $\liegood{h}(r,F')$ and
$[Y]_r=[Y']_r\in S_{\lie{g},\lie{h},r}(\f)$, then the fundamental
lemma (Equation~\ref{eqn:fl}) holds for $Y$ iff it holds for $Y'$.
\end{corollary}

\begin{proof}  The Equation~\ref{eqn:gfl} for the orbital
integrals depends on the local field $F$ only through the residue
field $\f$.
\end{proof}

\subsection{The ring of values for motivic integration}

Motivic integration takes values in a ring $\ring{K}$ defined by
Denef and Loeser \cite{ICM}. We briefly recall its definition, and
refer the reader to \cite{ICM} and \cite{Rat} for details.  First
of all, $K_0(\op{Var}_k)$ is the Grothendieck ring of varieties
over a field $k$ of characteristic zero.  It is generated by
symbols $[V]$, for every variety $V$ over $k$.  We omit the
relations. Let $\ring{L}=[\ring{A}^1_k]$ be the class of the
affine line. $K_0^{mot}(\op{Var}_k)$ is a quotient of
$K_0(\op{Var}_k)$ that is obtained by killing all
$\ring{L}$-torsion and by identifying $[V]$ and $[W]$ whenever
$[V]$ and $[W]$ are nonsingular projective varieties that become
equal in the category of Chow motives.  The ring $\mot$ is then
defined by inverting $\ring{L}$ and tensoring with $\ring{Q}$.  We
write $\ring{K}=\mot$ for a field $k$ that will be made explicit
below. We use $[\cdot]$ both for elements of the Grothendieck ring
$K_0(\op{Var}_k)$ and for their images in $\mot$.

Denef and Loeser also construct a completion of $\mot$.     This
completion is necessary in general, because integration is defined
as a limit. In the special setting that we consider, this
completion will not be necessary.  We will work exclusively with
the motivic volumes of {\it weakly stable}
subassignments.\footnote{Weakly stable subassignments will be
defined in Definition~\ref{def:weak}. There is an unfortunate
clash in terminology between `stable' in the sense of stable
conjugacy (stable orbital integrals, and so forth) and in the
sense of stable subassignments.  The context will make it clear
when stability is meant in the sense of weakly stable
subassignments.} As the completion will not be needed, we skip the
construction.

Let $k$ be the field of rational functions on
$S_{\lie{g},\lie{h},r}$. The generic fiber of each variety $U_i$
gives a element $[U_i]\in \mot$. It is natural to conjecture the
following geometric form of the fundamental lemma.  (Note that we
have cross-multiplied by the denominators in
Theorem~\ref{thm:variety}, to avoid a localization of the ring
$\mot$.)

\begin{conjecture}[Motivic fundamental lemma]\label{conj:gfl}
For each $(\lie{g},\lie{h})$ in Definition~\ref{def:pairs} and
each $r\in\ring{Q}$, we have the identity (in the ring
$\ring{K}$):
  \begin{equation}
  \begin{array}{lll}
  &{p(\lie{h},\lie{h},r)(\ring{L})}\sum_{i\in I(\lie{g},\lie{h},r)}
   b(\lie{g},\lie{h},r)_i
   [U(\lie{g},\lie{h},r)_{i}] =\\
   &{p(\lie{g},\lie{h},r)(\ring{L})} \sum_{i\in I(\lie{h},\lie{h},r)}
   b(\lie{h},\lie{h},r)_i
   [U(\lie{h},\lie{h},r)_{i}].
   \end{array}
  \end{equation}
\end{conjecture}

\subsection{Relation to the geometric fundamental
lemma}\label{sec:gfl}

\begin{remark}
Conjecture~\ref{conj:gfl} is closely related to the geometric
fundamental lemma described in \cite{GKM-fl} and \cite{Lfl}.
However, it is not clear whether our conjecture should be a
consequence of the geometric fundamental lemma as they formulate
it.  Their identity depends on $p$-adic parameters $\gamma$ and is
an identity built on varieties over the residue field $\f$.  Our
identity is a single universal identity (for each
$(\lie{g},\lie{h})$ and $r$) over the base field $k$, which is a
finitely generated extension of the field of rational numbers.
\end{remark}

\begin{remark}
The geometric approach to the fundamental lemma introduced by
Goresky, Kottwitz, and MacPherson assumes a local field of
positive characteristic.  Laumon has produced a proof of the
fundamental lemma for unitary groups (under a purity
hypothesis\footnote{The purity hypothesis has been verified for
{\it equivalued elements} in \cite{GKM-equiv}.}) over local fields
of positive characteristic \cite{Lfl}.  The unitary analogue of
Corollary~\ref{cor:F} (cf. Remark~\ref{remark:unitary}) extends
Laumon's proof of the fundamental lemma to local fields of
characteristic zero, at least for restricted elements.
\end{remark}

\begin{remark} In our geometric formulation of the fundamental lemma,
there is some loss of information when we pass from the ring
$\ring{Q}[S_{\lie{g},\lie{h},r}]$ to its field of fractions $k$.
It should be viewed as asserting that the fundamental lemma holds
generically.  We are forced to pass to the field of fractions
because the properties of the map $[\cdot]$ are based on the work
of Gillet and Soul\'e, which requires a field $k$ of
characteristic zero.
\end{remark}

\begin{remark} Although we have effective procedures for finding equations
for the data $M$, $I$, $b_i$, $U_i$, and $p$, it seems to be a
difficult problem in general to determine whether elements of the
ring $\ring{K}$ are equal. In particular, there is no known
decision procedure to determine whether the identity of
Conjecture~\ref{conj:gfl} is valid for a given $(\lie{g},\lie{h})$
and $r$.
\end{remark}

The rest of this paper is devoted to the proofs of the results
stated in this section.

\section{Characteristic Polynomials}

\subsection{$r$-reduction}

This section reviews some basic facts about polynomials and fields
extensions.
The proofs are elementary and are omitted.

Throughout the paper, we consider constants satisfying the
relations:
   \begin{equation}\begin{array}{lll}
   &r\in\ring{Q},\quad g,\ell,n,L,N\in\ring{Z};\quad N\ge1;\quad g\ge1;
   \quad r\ge0;\\
   &r=L/N;\ g=\gcd(L,N);\ \ell = L/g;\ n = N/g
   \end{array}
   \label{eqn:constants}
   \end{equation}

Let
\begin{equation}P = \lambda^N + \alpha_1 \lambda^{N-1} +\cdots + \alpha_n
\lambda^{N-n} + \cdots +
\alpha_{n g}\end{equation} be a polynomial with coefficients in
$F$ whose roots $\lambda_i$ in a fixed algebraic closure $\bar F$
satisfy
    \begin{equation}
    |\lambda_i| = q^{-r},
    \label{eqn:lambdaval}
    \end{equation}
for $i=1,\ldots,N$.  When Condition~\ref{eqn:lambdaval} holds, we
say that $P$ has {\it slope} $r$.

The coefficient $\alpha_{j}$ is a symmetric polynomial in
$\lambda_i$, which is homogeneous of degree $j$.  It follows from
Condition~\ref{eqn:lambdaval} that the coefficients of a
polynomial of {\it slope} $r$ satisfy
    \begin{equation}|\alpha_j|\le q^{-r j}.\end{equation}
In particular,
    \begin{equation}
    |\alpha_{n j}/\varpi^{\ell j}|\le 1.
    \end{equation}
Let $a_j$ be the image of the integer $\alpha_{n j}/\varpi^{\ell
j}$ in $\Fq$.

\begin{definition}
Set
   \begin{equation}
   R(\lambda) = \lambda^g + a_1 \lambda^{g-1}+\cdots+ a_g \in
   \Fq[\lambda],\quad
   \text{with } a_j = \op{res}_{\ell j} \alpha_{n j}.
   \label{eqn:r-reduct}
   \end{equation}
We call $R$ the {\it $r$-reduction\/} of $P$.  Let $t_r$ be the
map from $\{x\in \bar F \tq \op{ord}(x)=r\}$ to $\barf$ given by
   \begin{equation}
   t_r(\lambda) = ac (\lambda^n/\varpi^{\ell}) \in \barf^\times.
   \end{equation}
\end{definition}

\begin{lemma}
Let $P$ have slope $r$. If $\lambda$ is a root of $P$, then
    $t_r(\lambda)$ is a root of $R$.
\end{lemma}

\begin{definition}
We say that an integer in $\bar F$ is {\it topologically
unipotent\/} if its residue class is $1$.
\end{definition}

\begin{lemma}  Let $P$ have slope $r$.
Assume $p>n$, where $p$ is the characteristic  of $\Fq$. Assume
that $P$ has $N$ distinct roots $\lambda_1,\ldots,\lambda_N$.
Assume that $|\lambda_i-\lambda_j|=q^{-r}$ for all $i\ne j$. Then
the map $\lambda_j\mapsto t_r(\lambda_j)$ from roots of $P$ to
roots of $R$ is an $n$ to $1$ mapping onto the set of roots of
$R$.
\end{lemma}

\begin{corollary} The roots of $R$ are distinct.
\end{corollary}

We have a partial converse.

\begin{lemma}\label{lemma:converse}  Let $r$, $L$, $N$, $g$, $\ell$ be as
in Definition~\ref{eqn:constants}.  Assume that $p>n$, where $p$
is the characteristic of $\Fq$.  Let
    \begin{equation}P = \lambda^N + \alpha_1
\lambda^{N-1}+\cdots+\alpha_{ng}\end{equation}
be any polynomial in $F[x]$ such that $|\alpha_j|\le q^{-rj}$.
Define the $r$-reduction $R\in\Fq[\lambda]$ by
Condition~\ref{eqn:r-reduct}. Assume that $0$ is not a root of $R$
and that $R$ has distinct roots.  Then $P$ has slope $r$ and its
roots $\lambda_i$ satisfy $|\lambda_i-\lambda_j|=q^{-r}$ for all
$i\ne j$.
\end{lemma}

\begin{proof}
The inequality $|\alpha_j|\le q^{-rj}$ is strict when $n$ does not
divide $j$, because the left-hand side of the inequality is an
integral power of $q$.

Let
    \begin{equation}P_1 = \lambda^N + \alpha_{n} \lambda^{N-n}+\cdots
+\alpha_{ng}\end{equation}
be the polynomial obtained from $P$ by setting the coefficients
$\alpha_i$ to zero when $n$ does not divide $i$. Let
  \begin{equation}
  \tilde P(\lambda) = \varpi^{-r N}P(\lambda \varpi^r)\in \bar F[\lambda]
  \end{equation}
and similarly for $\tilde P_1(\lambda)$. The coefficients of
$\tilde P$ and $\tilde P_1$ are integers and the constant term is
congruent to $R(0)\ne0$.  It follows that the roots of $P$ (and
$P_1$) have absolute value $q^{-r}$. The resultant
    $\op{res}(\tilde P,\tilde P')$
is congruent modulo $\varpi$ to the resultant
    $\op{res}(\tilde P_1,\tilde P'_1)$.
Thus, the identity $|\lambda_i-\lambda_j|=q^{-r}$ follows if the
resultant for $\tilde P_1$ is a unit.  But $\tilde P_1$ has
nonzero roots and it is of the form $\dot R(x^n)$ where $\dot R$
is a lift to $F$ of $R$.  Hence the resultant for $\tilde P_1$ is
a unit iff the resultant $\op{res}(R,R')$ is nonzero.  This
follows from the assumption that $R$ has distinct roots.
\end{proof}

\subsection{Lifts of Polynomials}\label{sec:lift}

Let the constants $g,\ell,n,\ldots$ be related as in
Equations~\ref{eqn:constants}.   Let $R$ be a monic polynomial in
$\Fq[\lambda]$ of degree $g\ge 1$ with distinct nonzero roots:
    \begin{equation}R(\lambda) = \lambda^g + a_1 \lambda^{g-1} +\cdots+
a_g.\end{equation}

Let $\dot R$ be a lift to $F$.
    \begin{equation}\dot R(\lambda) = \lambda^g + \dot a_1\lambda^{g-1}
    +\cdots+ \dot a_g.\end{equation}
Thus, $\dot a_i$ is a representative in $\A$ of $a_i$ in $\f$.

\begin{lemma} If $R$ is irreducible, then $\dot R$ is irreducible.
\end{lemma}

\begin{proof} Gauss's lemma.
\end{proof}

Assume that $R$ is irreducible.  Let $F^{unr}_g$ be the unramified
extension of degree $g$ in $\bar F$.  The extension $F^{unr}_g$ is
a splitting field of $\dot R$ over $F$. Let $\dot\zeta$ be a root
of $\dot R$. Every root in $\bar F$ of the polynomial $x^n -
\varpi^{\ell}\dot\zeta$ generates a totally ramified extension of
degree $n$.  In particular, the polynomial is irreducible.
Consider the extension
    \begin{equation}F^{unr}_g((\varpi^\ell\dot\zeta)^{1/n})\iso
      F^{unr}_g[x]/(x^n-\varpi^\ell \dot \zeta).\end{equation}

Let $\dot \zeta_1,\ldots,\dot\zeta_g$ be the roots of $\dot R$.
The polynomial of degree $N$
    \begin{equation}\dot R_{(r)}(\lambda) = \prod_{i=1}^g (x^n -
\varpi^{\ell}\dot\zeta_i)
    \end{equation}
has coefficients in $F$.  The polynomial is irreducible over $F$.
The image of $\lambda^n/\varpi^\ell$ in the field extension
   \begin{equation}F[\lambda]/(\dot R_{(r)}(\lambda))\end{equation}
is a root of $\dot R$, which we use to identify $F^{unr}_g$ with a
subfield of this extension.

\begin{lemma}
\label{lemma:fieldiso}  For all
$\dot\zeta\in\{\dot\zeta_1,\ldots,\dot\zeta_g\}$,
   \begin{equation}
   F[\lambda]/(\dot R_{(r)}(\lambda)) \iso
   F^{unr}_g((\varpi^\ell\dot\zeta)^{1/n}).
   \end{equation}
In particular, as $\dot\zeta$ varies, the fields
\begin{equation}F^{unr}_g((\varpi^\ell\dot\zeta)^{1/n})\end{equation}
are isomorphic.
\end{lemma}

Let $\ddot R$ be another lift of $R$ to a degree $g$ monic
polynomial. Again, $F^{unr}_g$ is a splitting field of $\ddot R$.
Form $\ddot R_{(r)}(\lambda)$ as above.

\begin{lemma} \label{lemma:iso} Assume $p\not | \,n$.  The fields
    \begin{equation}
    F[\lambda]/(\ddot R_{(r)}(\lambda))\text{ and }
    F[\lambda]/(\dot R_{(r)}(\lambda))
    \end{equation} are isomorphic over $F^{unr}_g$. The fields
    \begin{equation}
    F^{unr}_g((\varpi^\ell\dot\zeta)^{1/n}) \text{ and }
    F^{unr}_g((\varpi^\ell\ddot\zeta)^{1/n})\end{equation} are isomorphic.
\end{lemma}

\begin{lemma} \label{lemma:dependsonly}
Assume $R$ is irreducible.  Let $P\in F[\lambda]$ be any monic
polynomial of degree $N$ with slope $r$ and with $r$-reduction
$R$.   Then $P$ is irreducible, and the isomorphism class of the
field extension
    \begin{equation}F[\lambda]/(P(\lambda))\end{equation}
depends only on $R$.
\end{lemma}

\begin{corollary} Let $R\in\f[\lambda]$ be an irreducible monic
polynomial of degree $g$.  Let $P$ have $r$-reduction  $R$. Assume
that $p$ is sufficiently large. The roots $\lambda_j$ of $P$ in an
algebraic closure satisfy
    \begin{equation}
    |\lambda_j| = q^{-r} \text{ and } |\lambda_i-\lambda_j|= q^{-r},
    \end{equation}
for $i\ne j$.
\end{corollary}

Now drop the assumption that $R$ is irreducible.  For each
irreducible factor $R_i$ of $R$ of degree $g_i$, the preceding
construction gives an unramified extension $F^{unr}_{g_i}$ of $F$
of degree $g_i$ and a totally ramified extension of
$F^{unr}_{g_i}$ of degree $n$. By Lemma~\ref{lemma:dependsonly},
these extensions of degree $n g_i$ are well-defined up to
isomorphism. Each factor $R_i$ has an ``$r$-lift'' $P_i$ of degree
$n g_i$. (Each $P_i$ is a monic polynomial with $r$-reduction
$R_i$.) Let $P$ be the product of the $P_i$. Its $r$-reduction is
$R$.

Start with a polynomial $P\in F[\lambda]$ that has slope $r$ and
with $r$-reduction $R$.  Since $P_i$ is irreducible if and only if
$R_i$ is, the factors of $P$ are of degree $n g_i$ and the factors
correspond in a $1$-$1$ fashion with the factors of $R$.

\subsection{Even Polynomials}\label{sec:even}

In this subsection, assume that $N$ is even and that
$P(-\lambda)=P(\lambda)$. That is, assume $P(\lambda) =
P^{(2)}(\lambda^2)$ for some polynomial $P^{(2)}$.  The constants
$g,\ell,n,\ldots$ continue to be defined as in
Definition~\ref{eqn:constants}. We show how to associate a
quadratic extension of algebras $F_i/F_i^\lb$ to each irreducible
factor of $P^{(2)}$.

If $n$ is also even then each pair $(\lambda,-\lambda)$ of roots
appear in the same fiber over the root $t_r(\lambda)$ of the
$r$-reduction $R$. For each irreducible factor of $R$, there are
totally ramified extension
\begin{equation}F_i^\lb=F^{unr}_{g_i}((\varpi^\ell\dot\zeta)^{2/n})\end{equation}
of degree
$n/2$ over $F^{unr}_{g_i}$. The extension $F_i/F_i^\lb$, where
   \begin{equation}F_i =F^{unr}_{g_i}((\varpi^\ell\dot\zeta)^{1/n}),
   \end{equation}
is a ramified quadratic extension.

If on the other hand, $n$ is odd, then each pair of roots
$(-\lambda,\lambda)$ is split between two fibers:\footnote{As
always, we assume that the residual characteristic is sufficiently
large, so that $p\ne2$.}
    \begin{equation}t_r(\lambda)\ne
-t_r(\lambda)=t_r(-\lambda).\end{equation}
The $r$-reduction $R$ is even. We write it as
$R(\lambda)=R^{(2)}(\lambda^2)$.  There are two types of
irreducible factors of $R$: those that are even polynomials and
those that are not. If $R_i$ is an irreducible factor that is an
even polynomial, then its degree $g_i$ is even.  A lift $R_i$ has
splitting field $F^{unr}_{g_i}$, with subfield $F^{unr}_{g_i/2}$.
The $(2r)$-reduction of $P^{(2)}_i$ is $R^{(2)}_i$. The field
extension
    \begin{equation}F_i^\lb = F[\lambda]/(P^{(2)}_i(\lambda))\end{equation}
can be identified with a
totally ramified extension of $F^{unr}_{g_i/2}$ of degree $n$. The
quadratic extension $F_i/F_i^\lb$, where
\begin{equation}F_i=F[\lambda]/(P_i(\lambda))\end{equation} is unramified.

If $n$ is odd and the irreducible factor $R_i$ is not an even
polynomial, then there is a matching irreducible factor $R_j$ with
$R_i(-\lambda)=R_j(\lambda)$.  The extension determined by $R_i$
in Lemma~\ref{lemma:fieldiso} is isomorphic to the extension
determined by $R_j$. We associate the algebra
   \begin{equation}
   F_i = F^{unr}_{g_i}((\varpi^\ell\dot\zeta)^{1/n})\oplus
   F^{unr}_{g_i}((\varpi^\ell(-\dot\zeta))^{1/n})\text{ over }
   F_i^\lb = F^{unr}_{g_i}((\varpi^\ell\dot\zeta)^{1/n})
   \end{equation}
with the factors $R_i$ and $R_j$.  The product $R_i R_j$ is of the
form $R^{(2)}_k(\lambda^2)$ for some irreducible factor
$R^{(2)}_k$ of $R^{(2)}$.

\section{Conjugacy in Classical Lie Algebras}\label{section: conjugacy}

The assumption remains in force that the characteristic of the
residue field is sufficiently large. (In particular, $p>2$, $p>n$,
and $p$ satisfies the restrictions of \cite{WA}.)

\subsection{Groups under consideration}

We consider symplectic and orthogonal groups.

In the symplectic case, we fix a nondegenerate skew form $q_V$ on
a vector space $V$ of even dimension $N$ over $F$.  We define
$Sp(q_V)$ to be the group preserving the form.  Concretely, we
assume that $V=F^N$, with $N$ even, and that $q_V$ is given by a
skew symmetric matrix $J$ on the standard basis $\{e_i\}$ of $F^N$
by
   \begin{equation}
     q_V(e_i,e_j)=\begin{cases}
      -1 & i+j=N+1,\ i>j\\
      1  & i+j=N+1,\ i<j\\
      0 & \text{otherwise}.
      \end{cases}
   \end{equation}
We let $\lie{sp}(N)$ be the corresponding lie algebra.

In the orthogonal case, we fix a nondegenerate symmetric form
$q_V$ on a vector space $V$ of dimension $d$ over $F$.  We let
$d=N$ for even orthogonal lie algebras and $d=N+1$ for odd
orthogonal lie algebras, where $N$ is even.  We define
$\lie{so}(d)$ to be the lie algebra associated with the form.  To
make things concrete, we take the vector spaces to be $F^d$ and
the symmetric forms to be defined by a matrix $J$ with respect to
the standard basis, where $J$ is the matrix given in \cite{WA} and
used in \cite{VTF}. We let $\lie{so}(d)$ be the corresponding lie
algebra. Its elements are $d\times d$ matrices that satisfy
   \begin{equation}
   {}^tX J + J X  = 0.
   \label{eqn:lie-relation}
   \end{equation}

\subsection{The parameter space $S_{\lie{g},r}$} \label{sec:parameter-space}

Assume that the constants $N$, $r$, $g,\ldots$ are related as in
Equation~\ref{eqn:constants}.  Assume $N$ is even.  Define
equivalence as in Definition~\ref{def:equiv}.  We define the
$r$-reduction $\lie{g}^{[r]}$ of a Lie algebra $\lie{g}$ in a
case-by-case manner.  It is defined in the following context. Let
\begin{equation}
  N = n g \text{ even};\ \lie{g} = \lie{sp}(N),\ \lie{so}(N),
\hbox{ or } \lie{so}(N+1).
\end{equation}
By Equation~\ref{eqn:lie-relation} we may take $\lie{g}$ to be
defined over $\ring{Z}$.  The $r$-reduction is again a lie algebra
over $\ring{Z}$, defined as follows:
\begin{equation}
\lie{sp}^{[r]}(N) = \begin{cases}
                \lie{sp}(g) & n \text{ odd}\\
                \lie{gl}(g) & n \text{ even}.
                \end{cases}
\end{equation}
\begin{equation}\lie{so}^{[r]}(N+1) = \begin{cases}
                \lie{so}(g+1) & n \text{ odd}\\
                \lie{gl}(g) & n \text{ even}.
                \end{cases}
\end{equation}
\begin{equation}\lie{so}^{[r]}(N) = \begin{cases}
                \lie{so}(g) & n \text{ odd}\\
                \lie{gl}(g) & n \text{ even}.
                \end{cases}
\end{equation}
The $r$-reductions are taken over $\ring{Z}$, but we also consider
them over $\ring{Q}$, $\f$, and so forth.

\begin{definition}
In the symplectic and odd orthogonal cases, we take
$S_{\lie{g},r}/\ring{Q}$ be the affine variety of regular
semi-simple conjugacy classes in $\lie{g}^{[r]}$; that is, the
adjoint quotient of the algebra $\lie{g}^{[r]}$.   In the even
orthogonal case when $n$ is odd, we take the subvariety of the
adjoint quotient of $\lie{g}^{[r]}$ on which the determinant (in
the standard representation of $\lie{g}$) is nonzero.  In the even
orthogonal case when $n$ are even, our construction is a bit more
exotic.  We take $S_{\lie{g},r}$ to be the affine variety of pairs
$(u,v)$ where $u$ is a regular semi-simple conjugacy class in
$\lie{gl}(g)$ with nonzero determinant, and $v$ is an element of
$\ring{G}_m$ such that $v^2 = -\det(u)$.
\end{definition}

\begin{example}  If $\lie{g}=\lie{sp}(6)$ and $r=1/3$, then
$n=3$, $\ell=1$, and $g=2$.  We have
   \begin{equation}
   \lie{g}^{[1/3]} = \lie{sp}(2) = \lie{sl}(2).
   \end{equation}
The set of regular elements are those with nonzero determinant:
   \begin{equation}
   \lie{sl}(2)\cap \op{GL}(2).
   \end{equation}
The conjugacy class (over $\bar{\ring{Q}}$) is determined by the
determinant.  The map
   $X\mapsto \det(X)$
induces an isomorphism $S_{\lie{g},r}\iso \ring{G}_m$ over
$\ring{Q}$.
\end{example}

\subsection{Pfaffians}\label{sec:pfaff}

In the case of even orthogonal Lie algebras, the stable conjugacy
class is not determined by the characteristic polynomial. Assume
that $P$ is the characteristic polynomial of a regular semi-simple
element.  Assume that $P$ has slope $r$.

We use the pfaffian to specify the map from the restricted  stable
conjugacy classes of slope $r$ of $\lie{so}(N)$ (over the $p$-adic
field) to stable conjugacy classes of $\lie{so}(g)$ (over the finite
field). The even orthogonal Lie algebra can be identified with $N$
by $N$ matrices satisfying
    \begin{equation}{}^t X J + J X= 0,\end{equation}
where $J$ is a symmetric matrix.  Then $JX$ is a skew symmetric
matrix.  Let $\op{pfaff}(JX)$ be its pfaffian.  (There is a
general discussion of pfaffians in \cite{G}.) The stable conjugacy
class is determined by the characteristic polynomial of $X$ and by
$\op{pfaff}(JX)$. We claim that
    $\det(J)= -1$.
In fact, the explicit choice of $J$  in \cite{WA} is a matrix with
$\pm1$ along the skew diagonal and zeroes elsewhere. Any symmetric
matrix of this form has determinant $-1$ (recall that $N$ is
even).

Let $x$ be an element in the $r$-reduction $\lie{so}^{[r]}(N)$
with characteristic polynomial $R$ and let $X$ be an element in
$\lie{so}(N)$ with characteristic polynomial $P$. Assume that the
$r$-reduction of $P$ is $R$.   Hence,
   \begin{equation}\det(X) = P(0) = \alpha_N,\quad \det(x) = R(0)=a_g =
   \op{ac}(\alpha_N).\label{eqn:det}\end{equation}

Assume that $n$ is odd. The $r$-reduction $\lie{g}^{[r]}$ is an
even orthogonal algebra $\lie{so}(g)$. In this case, let $j$ be
the symmetric matrix defining $\lie{so}(g)$. Assume it has the
same explicit form as $J$.  The matrix $j$ has determinant $-1$
for the same reasons as $J$. The square of the pfaffian is the
determinant.  It follows from Equation~\ref{eqn:det} that
    \begin{equation}
        \op{ac}\op{pfaff}(J X)^{2} = \op{pfaff}(j x)^{2}.
    \end{equation}
If $n$ is odd, then we take the matching conditions
$(X\leftrightarrow x)$ on stable conjugacy classes to be
    \begin{equation}
    \op{ac}\op{pfaff}(J X) = \op{pfaff}(j x).
      \label{eqn:pfaff1}
    \end{equation}

Assume that $n$ is even.  In this case, the $r$-reduction
$\lie{g}^{[r]}$ is the algebra
   \begin{equation}
   \lie{gl}(g).
   \end{equation}
But the variety $S_{\lie{g},r}$ consists of pairs
$x=(u,v)=(u(x),v(x))$ where $u$ is a regular semi-simple conjugacy
class in $\lie{gl}(g)$ with nonzero determinant and $v^2 =
-\det(u)$.  We have
\begin{equation}\op{ac}\op{pfaff}(J X)^2 = - a_g = -\det(u) =
v^2.\end{equation}
We take the matching condition to be
   \begin{equation}
   \op{ac}\op{pfaff}(J X) = v(x).
   \label{eqn:pfaff2}
   \end{equation}

\subsection{The map to $S_{\lie{g},r}$}

\begin{definition}  Let $f\in k[\lambda]$, where $k$ is any field.
Let $m$ be the multiplicity of the root $\lambda=0$ in $f$.  We
call the polynomial $f/\lambda^m$ the {\it nonzero part\/} of $f$.
\end{definition}

If $X$ is a regular semi-simple element in the symplectic lie
algebra, the nonzero part of its characteristic polynomial is the
same as the characteristic polynomial.  For odd orthogonal lie
algebras, the multiplicity is one, and for even orthogonal lie
algebras the multiplicity is zero or two.

Let $F$ be a $p$-adic field with residue field $\Fq$.  We
construct a map $\mu:\lie{g}(r) \to S_{\lie{g},r}(\f)$ when the
residual characteristic is sufficiently large. Let $X$ be a
restricted element of slope $r$. Let $P$ be the nonzero part of
the characteristic polynomial.  It has slope $r$.  Let $R$ be the
$r$-reduction of $P$.  The polynomial $R$ is the nonzero part of
the characteristic polynomial of an element in the reduced algebra
$\lie{g}^{[r]}$.  Its conjugacy class is an element of
$S_{\lie{g},r}$.  In the even orthogonal case, we add the
additional condition Equation~\ref{eqn:pfaff1} or
\ref{eqn:pfaff2}.  We use the same notation $\mu:\lie{h}(r)\to
S_{\lie{h},r}(\f)$ for the corresponding map for $\lie{h}$.

A stable conjugate of a restricted  element of slope $r$ is again
restricted  of slope $r$.  The map we have constructed depends
only on the stable conjugacy class of $X$, so that we may speak of
the image of a conjugacy class in $S_{\lie{g},r}(\f)$.  Two
elements $X$ and $X'$ are equivalent (in the sense of
Definition~\ref{eqn:alpha-equiv}) iff their stable conjugates are
equivalent. Thus, we may speak of equivalent stable conjugacy
classes.  If two stable conjugacy classes are equivalent, then
they define the same polynomial $R$ (and in the even orthogonal
case, the same pfaffian) and hence their images in $S_{\lie{g},r}$
are the same.

\begin{theorem} \label{thm:bij}
When the residual characteristic is sufficiently large,
the map $\mu$ induces a bijection between equivalence classes of
stable conjugacy classes of elements in $\lie{g}(r)$ and elements
of $S_{\lie{g},r}(\f)$.
\end{theorem}

\begin{proof}  We have already checked that $\mu$ induces a
well-defined map from equivalence classes of stable conjugacy
classes of elements in $\lie{g}(r)$ to $S_{\lie{g},r}(\f)$.

To see that this is onto, take the nonzero part $R$ of the
characteristic polynomial of $x\in S_{\lie{g},r}(\f)$.  Lift it to
an even polynomial $P=P^{(2)}(\lambda^2) = \dot R_{(r)}(\lambda)$
as in Section~\ref{sec:lift}. Associate with it a direct sum
$\oplus_{i\in I} F_i$ of algebras, as in Section~\ref{sec:even}.
These algebras embed as a Cartan subalgebra of $\lie{g}$ according
to the procedure given by Waldspurger \cite{WA}.  The polynomial
$P$ determines an element $X$ of this Cartan subalgebra such that
the nonzero part of its characteristic polynomial is $P$.  The
element $X$ belongs to $\lie{g}(r)$ and maps under $\mu$ to $x$.
(In the even orthogonal case, this is compatible with pfaffians).

To see that the map is $1$-$1$, we check that this lift from $x$
up to $X$ is well-defined up to stable conjugacy and equivalence.
The characteristic polynomial (together with the pfaffian in the
even orthogonal case) $P$ determines the stable conjugacy class in
these classical groups.  Also, the different lifts $\dot R_{(r)}$
give equivalent elements (Lemma~\ref{lemma:dependsonly}).  Thus
the theorem is established.
\end{proof}

\begin{corollary}\label{cor:conjugate}
If $X,X'\in\lie{g}(r)$ have the same image in $S_{\lie{g},r}(\f)$,
then their centralizers $G_X$ and $G_{X'}$ are stably conjugate.
\end{corollary}

\begin{proof}  The reduction $R$ determines the factorization of
$P_X$ and $P_{X'}$ and hence the direct sum of algebras
   $\oplus_{i\in I} F_i$ appearing in the proof of Theorem~\ref{thm:bij}.
According to \cite{WA}, this direct sum determines the stable
conjugacy class of the Cartan subalgebra (except in the even
orthogonal case, where the pfaffian must also be taken into
account).
\end{proof}

\subsection{Stable Orbital Integrals}\label{sec:stable}

Now comes the key result.  It allows us to parameterize stable
orbital integrals by the elements of $S_{\lie{g},r}(\f)$.
According to Corollary~\ref{cor:conjugate}, if $X$ and $X'$ have
the same image in $S_{\lie{g},r}(\f)$, then their centralizers are
stably conjugate. Hence we may normalize the orbital integrals of
$X$ and $X'$ by picking compatible measures on the centralizers of
$X$ and $X'$ (that is, we assume that conjugation from $G_X$ to
$G_{X'}$ preserves measures).

\begin{theorem} \label{thm:stable}
Let $X,X'\in\lie{g}(r)$.  Assume that $\mu(X)=\mu(X')$ in
$S_{\lie{g},r}(\f)$.  Assume that the measures on $G_X$ and
$G_{X'}$ are compatible in the sense just described.  Then the
stable orbital integral of $X$ is equal to the stable orbital
integral of $X'$.
\end{theorem}

\begin{proof}
By Corollary~\ref{cor:conjugate}, the Cartan subalgebras $G_X$ and
$G_X'$ are stably conjugate. Replacing $X'$ with a stable
conjugate, we may assume that $G_X = G_{X'}$.

Waldspurger parameterizes semi-simple elements, up to conjugacy,
by triples of data $(I,(a_i),(c_i))$ (up to an equivalence
relation) \cite{WA}.   Let
    $$
    (I,(a_i),(c_i)) \hbox{ and } (I',(a'_i),(c'_i))
    $$
be the parameters attached to $X$ and $X'$. Since $G_X = G_{X'}$,
we may assume that $I=I'$, $F_i' = F_i$, and $a_i,a'_i\in F_i$.
Since $X$ and $X'$ give the same element $\mu(X)=\mu(X')\in
S_{\lie{g},r}(\f)$, there is a unique bijection $\psi:I\to I$ and
isomorphisms $\rho_i:F_i\to F_{\psi(i)}$ such that $\rho_i(a_i)$
and $a'_{\psi i}$ have the same valuation and angular components
for each $i$. Passing to equivalent data, we may assume that
$\psi$ and $\rho_i$ are identity maps. Passing to a stable
conjugate of $X'$, we may assume that $c'_i=c_i$ for all $i\in I$.

After passing to this stable conjugate, we may assume that $X'\in
G_X$.  In fact, $G_X$ is determined by the data $(I,(a_i),(c_i))$
and $G_{X'}$ by the data $(I,(a'_i),(c_i))$.  According to the
criteria in \cite{WA}, we can take $G_X = G_{X'}$ if the element
$\epsilon(a_i,a'_i)=P'_X(a_i)/P'_{X'}(a'_i)$ is a norm of
$F_i/F_i^\lb$ for each $i$, where in general $P'_Y$ is the
derivative of the characteristic polynomial of $Y$.  However,
$X,X'\in\lie{g}(r)$ with the same image in $S_{\lie{g},r}$.
Moreover, each pair $(a_i,a'_i)$ has the same angular component
and valuation, . This implies that $\epsilon(a_i,a'_i)$ is
topologically unipotent in $F_i^\lb$ and therefore a square
whenever the residual characteristic is not $2$ (which we assume).
Thus we have $X'\in G_X$.

We may diagonalize the two elements of $G_X$ simultaneously and
compare the corresponding eigenvalues $\lambda_i$.  Since $a_i$
and $a_i'$ have the same valuation and angular components, there
exist topologically unipotent elements $u_i$ such that
    \begin{equation}\lambda_i(X) = u_i \lambda_i(X').\end{equation}

Pick an element $b$ in the building adapted to $X$ as in
Theorem~\ref{theorem: main}. According to Adler and Roche
\cite[\S~2]{A}, the Moy-Prasad filtration of a semi-simple Lie
algebra takes the following form on $G_X$:
    \begin{equation}
    \lie{g}_{b,r'} \cap G_X = \{ Y\in G_X \tq \forall\chi.\quad
    |\chi(Y)|\le q^{-r}\}
    \label{eqn:alpha}
    \end{equation}
and
    \begin{equation}\lie{g}_{b,r'+} \cap G_X = \{ Y\in G_X \tq
    \forall\chi.\quad
    |\chi(Y)|< q^{-r}\}.\end{equation}
Here $r'=c r$ for some positive scalar $c$ that translates between
the depth $r'$ and the slope $r$.  In these equations, $\chi$ runs
over all differentials of characters of the torus $T$ with Lie
algebra $G_X$.  In the concrete symplectic and orthogonal algebras
cases we consider, we can take $\chi$ to run over all integer
linear combinations of the nonzero eigenvalues $\lambda(Y)$ in the
standard representation.

Let $Z = X'-X$. If $\chi = \sum m_i\lambda_i$, then
    \begin{equation}\chi(Z) =
    \sum_i m_i(\lambda_i(X')-\lambda_i(X)) =
\sum m_i \lambda_i(X)(u_i-1).\end{equation} So $|\chi(Z)|<q^{-r}$.
It follows that $Z\in \lie{g}_{b,r'+}$.

We will show below in Equation~\ref{eqn:Dg} that the determinants
$D^{\gF,\lF}(X)$ and $D^{\gF,\lF}(X^\prime)$ in
Corollary~\ref{cor:unit element} are equal. Thus, by that
corollary, the orbital integrals of $X$ and $X' = X+Z$ are equal.

We can extend this result to stable orbital integrals as follows.
If $X$ corresponds to data $(I,(a_i),(c_i))$; and $X'$ corresponds
to the data $(I,(a'_i),(c_i))$.  Write $X_c$ and $X'_c$ to make
the dependence on the parameters $c=(c_i)$ explicit.  There is a
bijection between the orbits in the stable conjugacy class of $X$
and those in the stable conjugacy class of $X'$ given by
$X_c\leftrightarrow X'_c$.  Equation~\ref{eqn:alpha} gives slope
$r$ for $X_c$ and $X'_c$, which is independent of $c$.  The
argument given above for $X$ and $X'$ now applies for each $c$, to
give that the orbital integrals of $X_c$ and $X'_c$ are equal.
Summing over $c$, we find that the stable orbital integrals of $X$
and $X'$ are equal.
\end{proof}


\subsection{$\kappa$-orbital integrals}\label{section:kappa}

Let $(\lie{g},\lie{h})$ be one of the pairs of
Definition~\ref{def:pairs}.  Let $r\in \ring{Q}$.  The affine
variety $S_{\lie{h},r}$ is a product $S_{\lie{h_1},r}\times
S_{\lie{h_2},r}$, corresponding to the factors
$\lie{h}=\lie{h_1}\oplus \lie{h_2}$ of $\lie{h}$.

We have a map from stable conjugacy classes in
$\lie{h_1}\oplus\lie{h_2}$ to stable conjugacy classes in
$\lie{g}$ that is defined as follows.  Let
$Y=(Y_1,Y_2)\in\lie{h_1}\oplus\lie{h_2}$.  Let $P_1$ and $P_2$ be
the nonzero parts of the characteristic polynomials of $Y_1$ and
$Y_2$.  An image $X$ of $Y$ is an element whose characteristic
polynomial has nonzero part $P_1 P_2$.  In the even orthogonal
case, we also require that
    \begin{equation}
    \op{pfaff}(JX) = \op{pfaff}(JY_1)\op{pfaff}(JY_2)
    \end{equation}
(where each occurrence of $J$ is adapted to the appropriate size
of matrix).  $X$ is said to be an image of $Y$.

Let $S_{\lie{g},\lie{h},r}$ be the affine subvariety of
$S_{\lie{h},r}$ whose points are $\mu(Y)$ such that $Y$ has a
regular semi-simple image $X\in\lie{g}$.   The map $\mu$ restricts
to a map from the subset of $G$-regular elements of $\lie{h}(r)$
to $S_{\lie{g},\lie{h},r}$.

We have a morphism of varieties $S_{\lie{g},\lie{h},r}\to
S_{\lie{g},r}$ that is defined as follows.  If $y=(y_1,y_2)\in
S_{\lie{g},\lie{h},r}$, then there are corresponding nonzero parts
of characteristic polynomials $R_{y_1}$ and $R_{y_2}$.  The
element $y$ is mapped to $x\in S_{\lie{g},r}$ whose nonzero part
of the characteristic polynomial is $R_{y_1} R_{y_2}$. In the even
orthogonal case, we also add the condition that
   \begin{equation}\op{pfaff} (x j) = \op{pfaff} (y_1 j_1) \op{pfaff} (y_2
j_2)\end{equation}
or
   \begin{equation} v(x) = v(y_1) v(y_2)\end{equation}
as appropriate, where $v$ is the parameter of
Equation~\ref{eqn:pfaff2}.

\begin{definition}  We say that $X$ is an {\it image\/} of $y\in
S_{\lie{g},\lie{h},r}(\f)$ if $\mu(X)$ is the image of $y$ in
$S_{\lie{g},r}(\f)$.
\end{definition}

As in Section~\ref{sec:stable}, we pick compatible measures on $X$
and $X'$ when $G_X$ is stably conjugate to $G_{X'}$.

\begin{theorem}\label{thm:kappa}
For local fields of sufficiently large residual characteristic,
the $\kappa$-orbital integral $\orb{g}{h}(Y)$ of $Y\in \lie{h}(r)$
depends only on $\mu(Y)\in S_{\lie{g},\lie{h},r}(\f)$.
\end{theorem}

\begin{proof}  Consider $Y$ and $Y'$ that map to the
same parameter $y$.  The parameter $y$ determines $x\in
S_{\lie{g},r}(\f)$. Two element $X,X'\in\lie{g}(r)$ mapping to $x$
have stably conjugate centralizers $G_X$ and $G_{X'}$.  We assume
that $X$ is the image of $Y$ and that $X'$ is the image of $Y'$ in
$\lie{g}$.  Passing to a stable conjugate, we may assume (as in
the proof of Theorem~\ref{thm:stable}) that $G_X = G_{X'}$ and
that Waldspurger's parameters defining $X$ and $X'$ have the form
$(I,(a_i),(c_i))$ and $(I,(u_i a_i),(c_i))$ for some topologically
unipotent elements $u_i\in F_i$.  The constraints on the data
actually force $u_i\in F_i^\lb$.

The element $Y$ determines a partition of $I$ into two subsets and
the element $Y'$ determines a partition into two subsets. This
partition determines the $\kappa$.  We claim that the partition is
the same in both cases.  In fact, the partition is determined by
partitioning the even polynomial $P(\lambda) = P^{(2)}(\lambda^2)$
according to irreducible factors of $P^{(2)}$ (where $P$ as usual
is the nonzero part of the characteristic polynomial of $X$ or
$X'$).  The irreducible factors of $P^{(2)}$ are determined by the
irreducible factors of its $r$-reduction $R^{(2)}$ (see
Section~\ref{sec:even}), which is the same for both $X$ and $X'$,
since $X$ and $X'$ both map to $x$.  Thus, the partition is the
same in both cases.

We claim that for the chosen $X,X'$, we have
$\Delta(X,Y)=\Delta(X',Y')$, where $\Delta$ is the
Langlands-Shelstad transfer factor, as calculated by Waldspurger
in \cite[Ch.X]{WA}.  (In fact, if we index $X$ and $X'$ by the
data $(c_i)$ and write $X_c$, $X'_c$, we have $\Delta(X_c,Y) =
\Delta(X'_c,Y')$ for all parameters $c=(c_i)$.) For this, it
suffices to examine the explicit formula for the transfer factors
that Waldspurger calculates.  According to his calculations, the
ratio $\Delta(X_c,Y)/\Delta(X'_c,Y')$ is given as product of
characters of order $2$ on the following elements of $F_i^\lb$:
   \begin{equation}
   P'_X(a_i)/P'_{X'}(u_i a_i),
   \end{equation}
(where $P'$ is the derivative of $P$).  We claim that these elements
are topologically unipotent, so that the characters of order $2$ all
evaluate to $1$ on these elements. In fact, each $P'_X(a_i)$ is a
product of factors $\lambda_i(X)-a_i$ and $P'_X(u_i a_i)$ has the
corresponding form $\lambda_i(X')-u_i a_i$.  It follows from the
assumption that $X$ and $X'$ are restricted so that the quotient of
these two factors is topologically unipotent. Hence the claim.

It follows as in the proof of Theorem~\ref{thm:stable}, that the
orbital integrals of $X_c$ and $X'_c$ are equal for each
$c=(c_i)$.  Since the transfer factor is also the same for both
$X_c$ and $X'_c$, the $\kappa$-orbital integrals are equal.
\end{proof}

As a result of the theorem, henceforth we write $\orb{g}{h}(y)$
for the $\kappa$-orbital integral of the unit element for any
$Y\in\lie{h}(r)$ that maps to $y\in S_{\lie{g},\lie{h},r}(\f)$.

\section{The first order language of rings and Pas's language}
\label{sec:Pas}

The first order language of rings is a formal language in the
first order predicate calculus.  The concepts of logic and model
theory that we require in this paper can be found in Enderton
\cite{E} or Fried and Jarden, \cite{FJ}.

A language that is slightly more complicated language than the
first order language of rings is Pas's language. It is a
formalization of a fragment of the theory of Henselian fields. It
is described in \cite{Pas}, with additional comments in Denef and
Loeser's papers on motivic integration, particularly \cite{DLD},
and briefly in \cite{HC}. We briefly recall its most important
characteristics. The language is three-sorted in the sense of
\cite{E}. That is, quantifiers range over three distinct objects
that can be interpreted as a $p$-adic field $F$, its residue field
$\Fq$, and the additive group of values $\ring{Z}$ (or more
correctly, $\ring{Z}\cup\{+\infty\}$). The arithmetic of
$\ring{Z}$ is restricted to the additive theory. For the additive
theory of $\ring{Z}$, there is a procedure of quantifier
elimination due to Presburger \cite{Pres}. The language has
function symbols $\op{ac}$ and $\op{ord}$ that are interpreted in
a $p$-adic field as the angular component map and the valuation
function, respectively.

We recall the notion of a virtual set from \cite{VTF}.  It is a
syntactic extension of the first-order language.  Let $\calL$ be a
first order language (usually the first order language of rings,
Pas's language, or an extension of Pas's language obtained by
adjoining constants).  Let $\psi$ be a formula in $\calL$. We
write
    \begin{equation}`y\in \{x \tq \psi(x)\}\rq \text{ for }
`\psi(y).\rq\end{equation} The construct $\{x\tq \psi(x)\}$ is
called a {\it virtual set} (in the language $\calL$). Here, $x$ is
allowed to be a multi-variable symbol:
    $x = (x_1,\ldots,x_n)$, so that we have
    \begin{equation}
    `(y_1,\ldots,y_n)\in \{(x_1,\ldots,x_n) \tq \psi(x_1,\ldots,x_n)\}\rq
    \text { for }
    `\psi(y_1,\ldots,y_n)\rq
    \end{equation}
When we write $x\in {\bf A}$, it is to be understood that $x$ is a
vector of variable symbols, and that the length of that vector is
the number of free variables in the defining formula of ${\bf A}$.
Intersections, unions, complements and other standard operations
on sets can be applied to virtual sets.

For each of the split lie algebras considered in
Definition~\ref{def:pairs}, there is a virtual set (or {\it virtual
lie algebra\/}) defined by
   \begin{equation}\{X\tq {}^tX J + J X = 0\}\end{equation}
(viewed as a conjunction of equations in the free variables
$x_{ij}$).

\begin{lemma} \label{lemma:pas}
\begin{equation}\liegood{g}(r)\end{equation} is a virtual set in Pas's
language.
\end{lemma}

\begin{proof}
Let $P = \lambda^N + \alpha_1\lambda^{N-1} + \cdots + \alpha_{ng}$
be the nonzero part of the characteristic polynomial of $X$. For
each characteristic polynomial, the multiplicity of zero is known
in terms of the lie algebra $\lie{g}$.  Let $R = \lambda^g +
a_1\lambda^{g-1} + \cdots + a_g$ be the $r$-reduction of $P$. By
Definition~\ref{def:vg}, the restricted  elements are those such
that
   \begin{enumerate}
      \item $P$ is the nonzero part of the characteristic
      polynomial of $X$.
      \item The $r$-reduction of $P$ is $R$.
      \item $R$ has distinct roots.
      \item The multiplicity of $0$ in the characteristic
      polynomial of $X$ is $0$ or $1$.
      \item $|\alpha_{j}|\le q^{-rj}$ for all $j$.
      \item $a_g\ne 0$.
   \end{enumerate}
These conditions are all expressible in Pas's language.
\end{proof}

\subsection{Local Constancy}\label{sec:stability}

Let $\psi(x,\xi)$ be a formula in Pas's language, with free
variables $x=(x_1,\ldots,x_n)$ of the valued field sort and free
variables $\xi=(\xi_1,\ldots,\xi_m)$ of the residue field sort. We
set $|x|=n$ and $|\xi|=m$ to avoid a notational conflict with
Definition~\ref{eqn:constants}.

Given a formula $\psi(x,\xi)$, let $f_\psi$ be the auxiliary
formula
   \begin{equation}
   `\forall \xi\, x\, x'.\quad
   \op{ord} (x_i-x'_i)\ge M (\text{for }i=1,\ldots,n)
   \Rightarrow (\psi(x,\xi) \Leftrightarrow \psi(x',\xi))\rq.
   \end{equation}

\begin{definition}\label{def:loc-const}
We say that $\psi(x,\xi)$ is {\it locally constant} at level $M$
if $f_\psi(M)$ holds in all finite fields of sufficiently large
characteristic.
\end{definition}

\begin{remark}\label{remark:large-char}
Let us explain what it means for a formula $f_\psi(M)$ in Pas's
language to hold in all finite fields of sufficiently large
characteristic. When $M\in\ring{N}$ is substituted into $f_\psi$,
we obtain a sentence in Pas's language (with no free variables).
By quantifier elimination of the variables of the valued field and
value ring sort, we can find an equivalent sentence that contains
only terms of the residue field sort. (This involves discarding
finitely many primes as in \cite{HOI}.)  The formula $f_\psi$
holds in all finite field of sufficiently large characteristic, if
this equivalent sentence in the language of rings has a true
interpretation in all finite fields of sufficiently large
characteristic.
\end{remark}

\begin{lemma}\label{lemma:stable}
Suppose there exists $N$, such that if for every
$F$ whose residue characteristic is prime to $N$, there is an
$M_F$ (depending on $F$) such that the sentence $f_\psi^F(M_F)$
holds. Then $\psi(x,\xi)$ is stable of some level $M$.
\end{lemma}

\begin{proof}  Write
   \begin{equation}
   `f'_\psi(m)\rq \text{ for } `(\forall m' \ge m.\, f_\psi(m')) \wedge
      (\forall m' < m. \lnot f_\psi(m'))\rq.
   \end{equation}
It asserts that $m$ is the least level for which $\psi$ is stable
of level $M$.  When $F$ is such that there exists $M_F$ for which
$f^F_\psi(M_F)$ holds, then there is a unique $m_F$ for which
$f'_\psi(m_F)$ holds.  Let $\calF$ be the class of all $p$-adic
fields of sufficiently large residual characteristic.  Let
$\calF_N$ be the subset of these $p$-adic fields whose residue
characteristic is at least $N$.  By \cite[Thm.2]{HOI}, there
exists $N$ such that the set of natural numbers $\{m_F \tq
F\in\calF_N\}$ is bounded. Let $M$ be an upper bound. Then
$f^F_\psi(M)$ holds for all $p$-adic fields of sufficiently large
residual characteristic.
\end{proof}

This lemma can be immediately applied to the situation at hand.
For the classical groups we consider, the Langlands-Shelstad
transfer factor $\Delta$ is real valued.  We let
$\sign:\ring{R}\to \{-1,0,1\}$ be the usual sign function on the
Reals.   For each $\lie{g}$, $r\in\ring{Q}$, and
$\epsilon\in\{\pm1\}$, let $\psi_{\lie{g},r}^\epsilon$, be the
formula
   \begin{equation}
   \begin{array}{lll}
   `\psi_{\lie{g},r}^\epsilon(X,y)\rq \text{ for }&
   `X\in\lie{g}(r) \wedge \mu(X) = x \wedge y \in S_{\lie{g},r}
   \wedge\\
   &\quad \exists Y. \sign\Delta(X,Y) = \epsilon \wedge \mu(Y) = y \wedge
   (y\mapsto x)\rq.
   \end{array}
   \end{equation}
The expression $y\mapsto x$ indicates that $x$ is an image of $y$.

\begin{lemma}\label{lemma:psi-Pas}
The formula $\psi_{\lie{g},r}^\epsilon$ is expressible in Pas's
language.
\end{lemma}

\begin{proof}
The most difficult part of this claim is the assertion that
    $$
    \sign\Delta(X,Y)=\epsilon
    $$
is given by a formula in Pas's language. The main result of
\cite{VTF} shows that this part of the formula is actually given
in the first order language of rings. For the relation $y\in
S_{\lie{g},r}$ we use free variables of the residue sort,
constrained by the algebraic relations defining $S_{\lie{g},r}$ as
a subvariety of affine space.  For $\lie{g}(r)$, we use
Lemma~\ref{lemma:pas}.  The condition $X\in \lie{g}$ becomes
vacuous.  (The parameter $X$ is taken to be a set of $\dim\lie{g}$
variables ranging over $F$, under an identification of $\lie{g}$
with $F^{\dim\lie{g}}$.)  The result follows.
\end{proof}

\begin{definition} For $y\in S_{\lie{g},\lie{h},r}$, let
$\lie{g}(r)_y$ be the elements $X$ of $\lie{g}(r)$ such that $X$
maps to the image of $y$ in $S_{\lie{g},r}(\f)$.  We define
$\lie{h}(r)_y$ similarly.  Let $\lie{g}(r)^\epsilon_y$ be the
subset of $X\in \lie{g}(r)_y$ on which
$\sign\,\Delta(X,y)=\epsilon$.
\end{definition}

The formula $\psi_{\lie{g},r}^\epsilon(X,y)$ asserts that $X\in
\lie{g}(r)^\epsilon_y$.

\begin{lemma} \label{lemma:phi-stable}
$\psi_{\lie{g},r}^\epsilon$ is locally constant of some level $M$.
\end{lemma}

\begin{proof}  By Lemma~\ref{lemma:stable}, it is sufficient to show
that the interpretation of $\psi^\epsilon_{\lie{g},r}(\cdot,y)$,
for each $y\in S_{\lie{g},r}(\f)$ is locally constant as a
function of $X\in \lie{g}(F)$. The local constancy follows from
the local constancy (in $X$) of $\sign\Delta(X,y)$, which has
already been established (essentially) in the proof of
Theorem~\ref{thm:kappa}. In fact, $X$ corresponds to Waldspurger's
parameters $(I,(a_i),(c_i))$ and small perturbations of $X$ leave
$I$ and the fields $F_i$ unchanged.  Perturbing, $a_i\mapsto u_i
a_i$ where $u_i$ is topologically unipotent leaves the transfer
factor unchanged.  By the explicit formula for the transfer factor
in \cite{WA}, we have $\sign\Delta(X,y) = \sign\Delta(\op{Ad}\,g\,
X,y)$, and by Harish-Chandra's submersion principle such
conjugates of $G_X$ fill out a neighborhood of the regular element
$X$ in $\lie{g}$. The result follows.
\end{proof}

\section{Measures}

Each equivalence class of semi-simple orbits of slope $r$ forms an
open subset of the Lie algebra.  As a result, we may use the Lie
algebra form of the Weyl integration theorem to rewrite the
orbital integral as an integral over an open subset of the Lie
algebra (with the additive Haar measure).    We can express the
Lie algebra formulation of the fundamental lemma as an assertion
about volumes of regions in $\lie{g}(\A)$.

In the arguments that appear below, there are normalizations of
measures coming from three sources.  The first is the canonical
normalization of measures occurring in the theory of motivic
integration.  It is related to the Serre-Oesterl\'e measure that
arises in the integration theory of $p$-adic sets \cite{S},
\cite{O}. The second source of normalizations of measures comes
from the Weyl integration theorem.  The final source of
normalizations on measures comes from the fundamental lemma.  This
section shows that these various normalizations are compatible, in
the sense that the fundamental lemma takes on an appealing form
when the Weyl integration is used to express the fundamental lemma
as a statement involving the Serre-Oesterl\'e measures on the Lie
algebra.

\subsection{Normalization of Haar
measures}\label{sec:norm-measure}

Let $G$ be a reductive group with Lie algebra $\lie{g}$.  Assume
that $G$ and $\lie{g}$ are defined over $\A$. Normalize the
additive Haar measure $dX$ on $\lie{g}$ so that the lattice
$\lie{g}(\A)$ has volume $1$. Let $\varpi$ be a uniformizer of
$F$. The volume of the lattice $\varpi \lie{g}(\A)$ is
    \begin{equation}q^{-\dim(\lie{g})}.\end{equation}

Assume that the residue field characteristic is sufficiently
large.  Then there is a diffeomorphism between the lattice $\varpi
\lie{g}(\A)$ and a neighborhood $V$ of the origin in $G(\A)$
\cite{WH} as in Lemma~\ref{lemma: Adler 2}; see also \cite{deB}.
Explicitly, that neighborhood is the set of elements in $G(\A)$
with trivial image in $G(\bb{F}_q)$.  We normalize a Haar measure
$dg$ on $G$ so that the diffeomorphism preserves the volume of
$\varpi \lie{g}(\A)$ under the exponential map.

Thus,
    \begin{equation}
    \begin{array}{lll}
    \vol(G(\A),dg) &= [G(\A):V] \vol(\varpi  \lie{g}(\A),dX) \\
            &=
        {{|G(\bb{F}_q)|}{q^{-\dim\lie{g}}}}
    \end{array}.\end{equation}

\subsection{Serre-Oesterl\'e measures}

A general comparison theorem of Denef and Loeser \cite{DLD}, which
will be discussed further below, gives $p$-adic orbital integrals
as the trace of Frobenius on corresponding elements of the ring
$\ring{K}$. The normalization of $p$-adic orbital integrals in
their theorem is the canonical Serre-Oesterl\'e measure
\cite[Sec.8.2]{DLD}, \cite{S}, \cite{O}, \cite{V}.  Our
application will be to the integration of certain open subsets of
$\lie{g}(\A)$. By construction, the Serre-Oesterl\'e measure on
$\lie{g}(\A)$ is the additive Haar measure, normalized so that the
volume of $\lie{g}(\A)$ is $1$.

\subsection{The measures in the fundamental lemma}

Let $G$ be unramified, and let $K=G(\A)$ be a maximal compact
associated with a hyperspecial vertex.  The assertion of the
fundamental lemma requires a particular normalization of measures
\cite{HS}. Let $dg$ be a Haar measure on $G$ (such as the measure
given above), and let
    \begin{equation}f_{G} = {\frac{\chr
    G(\A)}{\vol(G(\A),dg)}}.
    \end{equation}
Descent implies that for $\gamma\in T(\A)$ a regular semi-simple
absolutely semi-simple element in a Cartan subgroup $T$, which is
regular modulo $\varpi$, we have \cite[Lemma 13.2]{HS}
    \begin{equation}\Phi^\kappa_{T,G}(\gamma,f_G) =
\frac{1}{\op{vol}(T(\A),dt)}.
    \end{equation}
The $\kappa$-orbital integral $\Phi^\kappa_{T,G}$ is computed with
respect to the quotient measure $dg/dt$. The $\kappa$-orbital
integral equals the stable orbital integral of the unit element on
the endoscopic group if the corresponding normalizations are used.

\subsection{Weyl integration formula}

Let $\calC(G)$ be a set of representatives for the conjugacy
classes of Cartan subgroups of $G$.  For each Cartan subgroup $T$,
let $\lie{t}$ be its Lie algebra.  Let $A_T$ be the split
component of $T$. Assume that $X\in\lie{g}$ has Jordan
decomposition
    \begin{equation}X = X_s + X_n,\end{equation} and let
    \begin{equation}D^{\gF}(X) = \det(\ad\,X |
\lie{g}/\lie{g}_{X_s}).\end{equation}

If $T$ is spit, normalize measures on $T$ to have volume $1$ on
the maximal compact subgroup of $T$.   Normalize measures on a
general Cartan $T$ by $\vol(T/A_T)=1$. Let the measures $dX$ and
$dg$ be compatibly normalized as in
Section~\ref{sec:norm-measure}. By \cite[page46]{WFT}, the Weyl
integration formula can be written as follows.
    \begin{equation}
    \begin{array}{lll}
    \int_\lie{g} f(X)dX &=
        \sum_{T\in\calC(G)}|W(G,T)|^{-1} \int_t
        |D^{\gF}(X)| \times\\
        &\int_{G/A_T} f(\Ad\,x X) dx\,dX.
    \end{array}
    \end{equation}

\subsection{Application}

An element $y\in S_{\lie{g},\lie{h},r}(\f)$ determines a Cartan
subalgebra $G_X$, up to stable conjugacy. For any such Cartan
subalgebra, we consider the volume
   \begin{equation}
   \frac{\vol(G_X\cap \lie{g}(r))}{|W(G,G_X)|}.
   \end{equation}
There is a corresponding volume in $\lie{h}(r)$.

\begin{lemma} \label{lemma:weyl}
For every $y\in S_{\lie{g},\lie{h},r}$, there is a constant
$\omega(y)$ such that for every $X\in\lie{g}(r)_y$ and every
$Y\in\lie{h}(r)_y$, we have
   \begin{equation}
   \omega(y) = \frac{\vol(G_X\cap \lie{g}(r))}{|W(G,G_X)|} =
   \frac{\vol(G_Y\cap \lie{h}(r))}{|W(H,G_Y)|}.
   \end{equation}
\end{lemma}

\begin{proof}  We prove the statement for $X\in\lie{g}(r)$.  The
proof for $Y\in\lie{h}(r)$ is similar and is left to the reader.
Fix a semi-simple element $X\in G_X$ that is an image of $y$. The
element $X$ has $|W(G,G_X)|$ conjugates in $G_X$. The set
$G_X\cap\lie{g}(r)_y$ is a disjoint union of $|W(G,G_X)|$ subsets
of $G_X$ indexed by the conjugates of $X$, consisting of elements
$\Omega_{X'}$ of $G_X\cap \lie{g}(r)$ closest to a given conjugate
$X'$. (That is, take Voronoi cells with centers at the conjugates
of $X$.) The volume of the set $\Omega_{X'}$ depends only on the
measure on the Cartan subalgebra $G_X$.  This volume is
independent of the ambient group. In particular, it is the same
for $G$ and $H$.
\end{proof}

For the classical lie algebras in this paper, there is a rational
number $a(r)$ depending on $r\in\ring{Q}$ such that the nonzero
values of the transfer factor on restricted  elements of slope $r$
have the form
   \begin{equation}
   \pm q^{a(r)}.
   \end{equation}
From Lemma~\ref{lemma:phi-stable}, we know that the sign $\pm$ is
given by a formula in the first order language of rings
\cite{VTF}.

On $\lie{g}(r)$ we have
    \begin{equation}\label{eqn:Dg}
    |D^{\gF}(X)| = |\prod_\alpha \alpha(X)| =
    q^{-r (\dim\lie{g}-\rank\lie{g})}.
    \end{equation}
Set $\delta_G = \dim\lie{g}-\rank\lie{g}$ and
   $\delta_H = \dim\lie{h}-\rank\lie{h}$.

By the proof of Theorem~\ref{thm:kappa}, the transfer factor at
$\Delta(X,Y)$ depends only on the image of $y$ in
$S_{\lie{g},\lie{h},r}(\f)$.  We write $\Delta(X,y)$ when
$Y\mapsto y$.  As above, let $\lie{g}(r)^\epsilon_y$ be the subset
of $\lie{g}(r)_y$ consisting of all $X$ such that the transfer
factor at $(X,y)$ is $\epsilon$.

\begin{lemma}\label{lemma:fld}
Let $(\lie{g},\lie{h})$ be one of the pairs of
Definition~\ref{def:pairs}.  Assume $r\ge 0$.  Let
    $$
    \langle G\rangle_q =|G(\ring{F}_q)| q^{-\op{dim}(G)},\quad
    \langle H\rangle_q =|H(\ring{F}_q)| q^{-\op{dim}(H)}.
    $$
    The fundamental lemma holds for $y\in S_{\lie{g},\lie{h},r}(\f)$ iff
    \begin{equation}
    \frac{\vol(\lie{g}(r)^+_y,dX)-\vol(\lie{g}(r)^-_y,dX)}
      {\langle G\rangle_q q^{-r\delta_G/2}} =
    \frac{\vol(\lie{h}(r)_y,dY)}
      {\langle H\rangle_q q^{-r\delta_H/2}}.
      \label{eqn:flq}
    \end{equation}
The measures are the Serre-Oesterl\'e measures $dX$ and $dY$ on
$\lie{g}(\A)$ and $\lie{h}(\A)$ respectively.
\end{lemma}

\begin{proof}  We use Lemma~\ref{lemma:weyl} to cancel the terms
$|W(G,T)|$ and $|W(H,T_H)|$ in the Weyl integration formula for
the algebras $\lie{g}$ and $\lie{h}$.   By the local constancy of
orbital integrals, the orbital integral equals the average of the
orbital integral over a neighborhood of the orbit inside
$\lie{g}(\A)$. By the Weyl integration formula this average is
equal to an integral in $\lie{g}(\A)$ with the additive Haar
measure on $\lie{g}(\A)$. By tracking the normalization of
measures, we obtain the result.
\end{proof}

\begin{example}  There is a simple case of Equation~\ref{eqn:flq} that can
be verified by hand.  This serves as a check on the correctness of
the normalization factors in the equation.  Assume $X$ is a
regular semi-simple element of $\lie{g}(\A)$ such that its image
in $\lie{g}(\f)$ is also regular semi-simple with $r=0$. A short
calculation shows that both sides in Equation~\ref{eqn:flq} equal
$1/|T(\f)|$, where $T$ is the centralizer of $X$ in $G$ (or $H$).
\end{example}

\section{Construction of Varieties}

We would like to apply the theory of motivic integration, as
developed in Denef and Loeser \cite{DLD} to Equation~\ref{eqn:flq}
to conclude the main result of the paper
(Theorem~\ref{thm:variety}). Unfortunately, the results of
\cite{DLD} do not give the desired results when there is a
parameterized family of integrals (in this case parameterized by
$y\in S_{\lie{g},\lie{h},r}(\f)$) rather than a single integral. The
forthcoming work of Cluckers and Loeser promises to give a general
theory of parameterized motivic integration \cite{CL}. However,
until those results become available, we confine ourselves to the
earlier papers of Denef and Loeser.

It is clear from an inspection of the proofs of \cite{DLD} that
the methods of that paper are not sufficient to show that general
parameterized families of integrals are ``motivic.''  However, if
we weaken the conclusions of their theorems slightly, the proofs
of that paper can be adapted to parameterized integrals.

Their paper must be adapted as follows.  Wherever they speak of an
element of the ring of motives $\ring{K}$, we speak instead of a
formal linear combination (with rational coefficients) of
varieties $U$ over $S = S_{\lie{g},\lie{h},r}$. Whenever they take
the trace of Frobenius on an element of $\ring{K}$, we count
points instead on the fiber $U_y$ over $y\in
S_{\lie{g},\lie{h},r}(\f)$.  With these slight modifications, we
can read through their proofs and check that the desired results
go through.  Note however, that they associate a canonically
determined element of $\ring{K}$ to definable subassignments, but
our representation as a linear combination of varieties is far
from unique.

We give a few technical details in the paragraphs that follow
about how specific arguments in their paper are to be adapted to
the parameterized orbital integrals in this paper.

We make use of the following variant of one of the main results of
\cite{DLD}.

\begin{theorem}\label{thm:existU}  Let $\psi(x,\xi) =
\psi(x_1,\ldots,x_n,\xi_1,\ldots,\xi_m)$ be a formula in Pas's
language with free variables $x=(x_1,\ldots,x_n)$ of the
valued-field sort, free variables $\xi = (\xi_1,\ldots,\xi_m)$ of
the residue field sort, and no other free variables.  Assume that
$\psi(x,\xi)$ is locally constant of some level.  Let $S\subset
\ring{A}^m_{/\ring{Z}[1/M]}$ be an affine variety.  Assume that
$\psi$ projects to $S$ in the sense that the following sentence in
the first order language of rings holds for all finite fields of
sufficiently large characteristic\footnote{This is meant in the
sense of Remark~\ref{remark:large-char}.} (in particular
$(q,\ell)=1$):
   \begin{equation}
   \forall\xi.\quad
   (\exists x.\, \psi{(x,\xi)})
   \Rightarrow (\xi\in S)
   \label{eqn:S}
   \end{equation}
Then there exist a natural number $M$ (with $\ell | M$), a finite
indexing set $I$, constants $b_i\in\ring{Q}$ for $i\in I$,
varieties $U_i$ over $S$, and a polynomial $p(x)\in\ring{Q}[x]$ of
the form
   \begin{equation}p(x) = x^\ell \prod_{i=1}^{\ell'}
    (x^{k_i}-1)
    \end{equation} with the following property.
\begin{itemize}
   \item
   For all $p$-adic fields $F$ and all residue fields $\f$, such
   that $(q,M)=1$, and for all $y \in S(\f)$, we have
   \begin{equation}
   \op{vol}\,(\{ x\in O_F^n \tq \psi^F(x,y)\},dx) =
   \frac{1}{p(q)}\sum_{i\in I} |U_{i,y}(\f)|,
   \label{eqn:volU}
   \end{equation}
where $dx$ is the additive Haar measure on $O_F^n$ normalized so
that $\op{vol}\,(O_F^n) = 1.$
\end{itemize}
\end{theorem}

\begin{corollary} Under the same hypotheses, for all $p$-adic
fields $F$, all residue fields $\f$ such
   that $(q,M)=1$, and for all $y \in S(\f)$,
    $$
    \op{vol}\,(\{ x\in O_F^n \tq \psi^F(x,y)\},dx)
    $$
    depends on $F$ only through $\f$.
\end{corollary}

We supply a sketch of the proof of the theorem, with references to
\cite{DLD} for details.  Our argument relies on many of the ideas
and constructions from \cite{DLD}.  The rest of this paper follows
that paper closely; and our argument should be read with that
paper at hand. Before turning to the proof, we give several
reductions.

\subsection{Reduction to covers} First we note that there is no
loss of generality in working with coordinate patches that cover a
variety.  In fact, if $X_{/S}$ is any variety with cover
$\{U_j\}_{j\in I}$ (with each $U_i\to S$), then
   \begin{equation}
   \forall y\in S(\f).\quad
   |X_y(\f)| = \sum_{J\subset I} (-1)^{|J|} |\cap_{j\in J} U_{j,y}(\f)|.
   \end{equation}
This equation can be used to combine the results obtained on
separate coordinate patches.

\subsection{Reduction to weakly stable subassignments} If
$\psi(x,\xi)$ is any formula, let $\psi_N(x,\xi)$ be the formula
    \begin{equation}
    \exists x'.\quad \op{ord}(x_i-x_i')\ge N (\text{for } i=1,\ldots,n)
    \land \psi(x',\xi).
    \end{equation}
The formula $\psi_N$ is true at $x$ whenever $x$ is `close' to
$x'$ that satisfies $\psi$.  Assume that $\psi$ is locally
constant of some level $N$.  Then $\psi_N$ is also locally
constant of level $N$.  If $\psi$ satisfies Equation~\ref{eqn:S},
then $\psi_N$ does too, because
    \begin{equation}
    \forall \xi.\ [
    (\exists x.~
    \psi(x,\xi))~~\Leftrightarrow~~
    (\exists x.~\psi_N(x,\xi))].
    \end{equation}
Moreover, for all $p$-adic fields $F$ of sufficiently large
residue characteristic, $\psi^F(x,\xi)\Leftrightarrow
\psi_N^F(x,\xi)$.  Thus, it is enough to prove
Theorem~\ref{thm:existU} for $\psi_N$ rather than $\psi$.

\begin{definition}
Let $h:C\to\op{Set}$ be a functor into the category of sets.  A
{\it subassignment\/} of $h$ is a function $f$ from the objects of
$C$ such that $f(C)\subset h(C)$.
\end{definition}

Let $C$ be the category of fields of characteristic zero.  For a
fixed $m,n$, let $h=h_{m,n}:C\to\op{Set}$ be given by
    $$k\mapsto k[[t]]^m\times k^n.$$
A formula $\psi(x,\xi)$ defines a subassignment $f$ of $h_{m,n}$
(with $(m,n)=(|x|,|\xi|)$ by
    $$
    f(k) = \{(x,\xi)\in k[[t]]^m\times k^n\tq \psi^k(x,\xi)\}.
    $$
\begin{definition}\label{def:weak}
A subassignment is {\it definable\/} if it is attached to a
formula in this way.   A definable subassignment (attached to
$\psi$) is {\it weakly stable\/} of level $N$ if for every field
$k$ of characteristic zero
    \begin{equation}
    \forall \xi\in k^n.~
    \forall x~x'.~
    (\forall i.~\op{ord}_k(x_i-x_i')\ge N) ~\Rightarrow~
    [\psi^k(x,\xi) \Leftrightarrow \psi^k(x',\xi)].
    \end{equation}
\end{definition}

This is well-defined: if a given subassignment is attached to both
$\psi$ and $\psi'$, it is weakly stable for $\psi$ iff it is
weakly stable for $\psi'$.  For any $\psi$ in Pas's language, the
subassignment of $\psi_N$ is weakly stable of level $N$. Thus, we
reduce to the case where $\psi_N$ determines a weakly stable
subassignment.

\subsection{Reduction to special formulas}

\begin{definition}
    $\psi(x,\xi)$ is a {\it special formula of bounded
    representation\/}
    if it can be expressed as a boolean combination of formulas
\begin{equation}
   \theta(\xi,\op{ac}\,f_1(x),\ldots,\op{ac}\,f_{m'}(x))
   \land (\op{ord} f_1(x) = N_1)
   \land \cdots
   \land (\op{ord} f_{m'}(x) = N_{m'})
   \label{eqn:special}
   \end{equation}
where each $N_i\ne 0$ and $\theta$ is a formula in theory of rings
in the variables and constants of the residue field sort.
\end{definition}

It is easy to see that each special formula of bounded
representation determines a weakly stable subassignment.  (This
follows from the fact that the functions $\op{ord}(f(x))$ and
$\op{ac}(f(x))$ are locally constant in $x$ when $\op{ord}(f(x))$
is a fixed, nonzero integer.)  We have the following converse.

\begin{lemma}  Every weakly stable subassignment is defined by a
special formula of bounded representation.
\end{lemma}

\begin{proof}  Assume that the weakly stable subassignment is
defined by a formula $\psi(x,\xi)$.  Apply quantifier elimination
(following Pas and Presburger) to eliminate all quantifiers of the
valued field sort and the value group sort in the formula
$\psi(x,\xi)$.  In Presburger quantifier elimination, the additive
language of the integers is augmented by function symbols for
congruences modulo $n$ for each $n$.  Each formula can be written
in disjunctive normal form.  Each disjunct is a conjunct of three
formulas one for the valued field sort, one for the value group
sort, and one for the residue field sort.   The conjunct for the
valued field sort can be eliminated, for example, by replacing
$f(x)=0$ with $\op{ac}(f(x))=0$ (as an extra condition in the
conjunct of residue field sort).  What results is a so-called
special formula; that is, a formula that can be expressed as a
boolean combination of formulas
    \begin{equation}
    \begin{array}{lll}
   &\op{ord} f_1(x) \ge \op{ord} f_2 (x) + a,\\
   &\op{ord} f_1(x) \equiv a \, \mod b,\\
   &\theta(\xi,\op{ac}\,f_1(x),\ldots,\op{ac}\,f_{m'}(x))
   \end{array}
    \end{equation}
If we show that each $f_i$ that appears can be assumed to satisfy
a bound $\op{ord}(f_i(x)) < N_i$ for some $N_i\ne0$, then the
lemma follows, by breaking each special formula into a disjunction
of finitely many cases, according to the possible values of
$\op{ord}(f_i)$.  This result is essentially identical to a lemma
of Denef and Loeser \cite[Lemma~2.8]{DLG}. (Denef and Loeser
assume that $\theta$ contains no free variables $\xi$. However, it
is trivial to check that their proof goes through without
modification in this slightly more general setting.)
\end{proof}

Not only can we reduce to special formulas of bounded
representation, but we can also reduce to a single formula like
Formula~\ref{eqn:special}.  (That is, no boolean combinations are
required.)  First of all, if $B(\psi_1,\ldots,\psi_k)$ is a given
boolean polynomial, then we obtain the same subassignment if we
replace each $\psi_i(x,\xi)$ with $\psi_i(x,\xi)\land
\phi_S(\xi)$, where $\phi_S(\xi)$ is the formula that asserts that
$\xi\in S$. Thus, there is no loss in generality in taking the
boolean operations ``relative to $S$''.  Consider conjunction.  A
conjunction of formulas of the form (Formula~\ref{eqn:special}) is
again of the same form. If we have a disjunction
$\psi_1\lor\psi_2$ of this form, and if we can prove
Theorem~\ref{thm:existU} for $\psi_1$, $\psi_2$,
$\psi_1\land\psi_2$ (with $A_1$, $A_2$, and $A_{12}$ as the
right-hand side of Equation~\ref{eqn:volU}), then we have
Theorem~\ref{thm:existU} for $\psi_1\lor\psi_2$ (with
$A_1+A_2-A_{12}$ as the right-hand side of
Equation~\ref{eqn:volU}). Finally, if we have the negation
$\lnot\psi(x,\xi)$ (relative to $S$) of a special formula, we use
    $$
    (\lnot\psi\land\phi_S)\lor(\psi\land\phi_S) = \phi_S.
    $$
to eliminate $\lnot\psi$.

Now we are ready to move to the proof of the representation
theorem for formulas.

\begin{proof} (Theorem~\ref{thm:existU}).
Let $F = \prod f_i$, the product extending over the functions
$f_i$ appearing in the representation of $\psi$ in
Formula~\ref{eqn:special}.   As in
\cite[Proof~7.1.1,Proof~8.3.1]{DLD}, take an embedded resolution
of $F=0$.  This resolution is independent of $\xi$.  It is {\it
good\/} in the sense of \cite{D-degree} when the residual
characteristic is sufficiently large.

The function $F$ comes from an expression in the first order
theory of rings.  We may thus interpret it as a polynomial in
$\ring{Q}[x]$, and the embedded resolution as a resolution of
$F=0$ over $\ring{Q}$.  In each suitably chosen coordinate patch
$W$ in the resolution, $F=0$ defines a divisor with normal
crossings. If $F$ is given in local coordinates as $\alpha
u_1^{k_1}\cdots u_n^{k_n}$, we let for $I \subset \{1,\ldots,n\}$
   \begin{equation}
   E_I = \{(u_i)\in W \tq  u_i = 0 \Leftrightarrow i\in I\}.
   \end{equation}
On each coordinate patch $W$, for each $I$, we obtain a formula
$\theta_I$ in the first order language of rings as follows.  Pull
back each $f_j$ to a function $w_j$ on the resolution.  If $w_j$
is identically zero on $E_I$, let $w'_j = 0$; otherwise, let $w'_j
= w_j$. Set
   \begin{equation}\theta_I(w,\xi) = (w \in E_I) \wedge \theta(w',\xi),
   \end{equation}
where $\theta$ is formula in the first order theory of rings
defining the special formula of bounded representation.

\smallskip
We construct a Galois stratification for the formula $\theta_I$ as
in \cite{DLD}.  (See \cite{FJ} for a review of Galois
stratifications.)  The particular version of Galois stratification
that we use is that of Lemma~\ref{lemma:FJ} below. The varieties
$U_i$ are constructed from individual strata of the Galois
stratifications of $\theta_I$, then summing over all strata for
all $I$.

Let $(C/A,\op{Con})$ be a colored Galois cover with Galois group
$G$ that arises in the Galois stratification of some $\theta_I$.
We assume that $C$ and $A$ are affine. By Artin induction, the
central function of $G$ given by
   \begin{equation}\alpha(x) = \begin{cases} 1 & \text{if } \langle x
\rangle
   \in \op{Con}\\
         0 & \text{otherwise.}
         \end{cases}
   \end{equation}
is a rational linear combination $\alpha = \sum n_H \op{Ind}_H^G
1_H$ of characters induced from trivial characters on cyclic
subgroups $H$ of $G$.  The formal linear combination of varieties
that corresponds to this colored Galois cover is
   \begin{equation}\sum n_H [C/H].\end{equation}

The morphism $U_i \to S$ is the composite
   \begin{equation}C/H \to C/G = A \to \ring{A}^{|u|+|\xi|} \to
   \ring{A}^{|\xi|}.\end{equation}
The last morphism is projection onto the last $|\xi|$ factors.
Recall that $A$ is a stratum in some $E_I$, which is a
constructible subset of a coordinate patch $W$ in the resolution.
We use the coordinate functions $u_i$ of the coordinate patch as
the free variables in the formula.  We may take $S$ be given as a
affine variety in $\ring{A}^{|\xi|}$.  The image of $U_i$ in
$\ring{A}^{|\xi|}$ lies in $S$ because of the
assumption~\ref{eqn:S} in the statement of the theorem.

For each $(u,\xi)\in \ring{A}^{|u|+|\xi|}(\f)$, we have by
\cite[3.3.2]{DLD}
   \begin{equation}
   \sum n_H |(C/H)_{u,\xi}(\f)| = \begin{cases}
      1 & \text{if } \theta_I(u,\xi)\\
      0 &\text{otherwise.}
      \end{cases}
   \end{equation}
Summing over $u$, we find that
   \begin{equation}
   \forall\xi.\quad
      \sum n_H |(C/H)_\xi(\f)| = |\{u \tq \theta_I(u,\xi)\}|.
   \end{equation}
This expresses the number of solutions $u\in\f^{|u|}$ of the
formula $\theta_I$ for a given $\xi$ as a linear combination of
number of points on varieties, with varieties that are independent
of the element $\xi$.  This is the essential point of the proof.
The rest of the argument is no different from that \cite{DLD}.
\end{proof}

The following result, which was used in the proof of
Theorem~\ref{thm:existU}, is taken directly from Fried and Jarden.

\begin{lemma} \label{lemma:FJ} (\cite[Prop.~26.7]{FJ})
Let $k$ be a nonzero integer and let $\theta(Y)$ be the Galois
formula
   \begin{equation}
   (Q_1 X_1)\ldots (Q_m X_m) [\op{Ar}(X,Y) \subseteq
   \op{Con}(\calA(\calC))]
   \end{equation}
with respect to a Galois stratification $\calA(\calC)$ of
$\ring{A}^{m+n}$ over $\ring{Z}[1/k]$, where $\calC$ is the family
of all finite cyclic groups.  Then we can effectively find a
nonzero multiple $\ell$ of $k$, and a Galois stratification
$\calB(\calC)$ of $\ring{A}^n$ over $\ring{Z}[1/\ell]$ such that
the following holds: for each finite field $M$ of characteristic
not dividing $\ell$, and for each $b\in \ring{A}^n(M)$,
   \begin{equation}
   M \models \, \theta(b) \Leftrightarrow \op{Ar}(\calB,M,b)\subseteq
   \op{Con}(\calB(\calC)).
   \end{equation}
\end{lemma}

\begin{remark}\label{remark:unitary}
We have accomplished our objective of showing that the orbital
integrals of restricted elements in symplectic and orthogonal
algebras count points on varieties over finite fields. There
should be no significant difficulties in extending these results
to non-split cases such as the unitary algebra
$\lie{g}=\lie{u}(c)$ and $\lie{h}=\lie{u}(a)\oplus\lie{u}(b)$. The
necessary analysis of the Langlands-Shelstad transfer factor has
already been carried out in \cite{VTF}.  That paper describes how
to extend Pas's language by a function symbol whose interpretation
is conjugation with respect to a separable quadratic extension
$E/F$. It describes how to expand formulas in the extended
language into formulas in Pas's language.  This trick can be used
to treat algebras that split over a quadratic extension.
\end{remark}

\section{Appendix : local constancy of orbital integrals as a
Corollary of \cite{KM}[3.1.9], by Julee Kim}

In general, we will keep the notation from \cite{KM}.

We assume the residue characterstic $p$ of $k$ is sufficiently
large (for more precise condition, we refer to \cite{KM}[1.4]).
Let $\Gamma$ be a good semisimple element of depth $r$, and let
$\mathbf G'=\mathbf C_{\mathbf G}(\Gamma)$. Let $\mathcal
O(\Gamma)$ be the set of $G$-orbits in $\mathfrak g$ whose
closures contain $\Gamma$. Denote the subspace $\sum_{x\in\mathcal
B(\mathbf G,k)} C_c(\mathfrak g/\mathfrak g_{x,r^+})$ of
$C_c^\infty(\mathfrak g)$ by $D_{r^+}$.

Recall the following theorem in \cite{KM}[3.1.9]:

\begin{theorem}
Let $\mathcal J(\Gamma+\mathfrak g'_{r^+})$ be the set of
$G$-invariant distributions supported on $^G(\Gamma+\mathfrak
g'_{r^+})$ and let $\mathcal J_\Gamma$ be the span of orbital
integrals associated to orbits in $\mathcal O(\Gamma)$. Then
\[
\mathcal J(\Gamma+\mathfrak g'_{r^+})\equiv\mathcal J_\Gamma
\]
on $D_{r^+}$.
\end{theorem}

From now on, we also assume that $\Gamma$ is {\bf either regular
or elliptic}. Let $x\in\mathcal B(\mathbf G',k)$, and let
$\Gamma'\in\Gamma+\mathfrak g_{x,r^+}$ be a good element. Since we
have $\Gamma+\mathfrak g_{x,r^+}=\,^{G_{x,0^+}}(\Gamma+\mathfrak
g'_{x,r^+})$ by \cite{KM}[2.3.5], for the purpose of comparing the
orbital integrals associated to $\Gamma$ and $\Gamma'$, we may
assume $\Gamma'\in\Gamma+\mathfrak g'_{x,r^+}$. Then it follows
from \cite{KM}[2.3.6] that
\[
G'=C_G(\Gamma)=C_G(\Gamma').
\]


Fix a Haar measure on $G/G'$, and denote the orbital integrals
associated to $\Gamma$ and $\Gamma'$ by $\mu_\Gamma$ and
$\mu_{\Gamma'}$ respectively.

\begin{theorem}
Let $\Gamma$ and $\Gamma'$ be as above. Then
$\mu_\Gamma(f)=\mu_{\Gamma'}(f)$ for any $f\in D_{r^+}$.
\end{theorem}

\proof Note that $\mu_{\Gamma'}\in \mathcal J(\Gamma+\mathfrak
g'_{r^+})$. Since $\mathcal O(\Gamma)$ has a single element
$^G\Gamma$, by the above theorem, $\mu_{\Gamma'}\equiv
c\cdot\mu_\Gamma$ on $D_{r^+}$ for some constant $c$.

Let $f_\Gamma\in D_{r^+}$ be the characteristic function supported
on $\Gamma+\mathfrak g_{x,r^+}=\Gamma'+\mathfrak g_{x,r^+}$. Then,
we have $\mu_\Gamma(f_\Gamma)=\,vol_{G/G'}(G_{x,0^+}G'/G')
=\mu_{\Gamma'}(f_\Gamma)$. Hence, $c=1$ and
$\mu_\Gamma=\mu_{\Gamma'}$ on $D_{r^+}$. \qed

\bibliographystyle{amsplain}

    \renewcommand{\thefootnote}{}
    \footnote{The research of the second author was supported in part by the
NSF.}

    \footnote{Copyright 2003, Clifton Cunningham and Thomas C. Hales.
This work is licensed under the Creative Commons Attribution
License. To view a copy of this license, visit
http://creativecommons.org/licenses/by/1.0/ or send a letter to
Creative Commons, 559 Nathan Abbott Way, Stanford, California
94305, USA.}

\end{document}